\documentclass[12pt,a4paper]{amsart}
\usepackage[a4paper,margin=1in]{geometry}
\usepackage{amsmath,amssymb,amsthm,mathtools}
\usepackage[protrusion=true,expansion=false]{microtype}
\usepackage{enumitem,booktabs}
\usepackage[colorlinks=true,linkcolor=blue!50!black,citecolor=blue!50!black,urlcolor=blue!50!black]{hyperref}
\usepackage{xcolor}

\setlist{itemsep=3pt,topsep=5pt}
\allowdisplaybreaks

\theoremstyle{plain}
\newtheorem{mainthm}{Theorem}

\newtheorem*{thmnn}{Theorem}
\newtheorem{theorem}{Theorem}[section]
\newtheorem{lemma}[theorem]{Lemma}
\newtheorem{proposition}[theorem]{Proposition}
\newtheorem{corollary}[theorem]{Corollary}

\theoremstyle{definition}
\newtheorem{definition}[theorem]{Definition}
\newtheorem{example}[theorem]{Example}
\theoremstyle{remark}
\newtheorem{remark}[theorem]{Remark}

\newcommand{\QQ}{\mathbb{Q}}
\newcommand{\ZZ}{\mathbb{Z}}
\newcommand{\CC}{\mathbb{C}}
\newcommand{\FF}{\mathbb{F}}

\newcommand{\Qbar}{\overline{\mathbb{Q}}}
\newcommand{\Kbar}{\overline{K}}
\newcommand{\hh}{\widehat{h}}
\newcommand{\OO}{\mathcal{O}}
\newcommand{\OK}{\mathcal{O}_K}
\newcommand{\p}{\mathfrak{p}}
\newcommand{\Sfin}{\mathcal S^{\mathrm{fin}}_\p}
\DeclareMathOperator{\Res}{Res}
\DeclareMathOperator{\disc}{disc}
\DeclareMathOperator{\Nm}{N}
\begin{document}

\title[Polynomials of Frobenius type]{Dobrowolski-type lower bounds for
canonical heights of polynomials of Frobenius type}
\author{Jorge Mello}
\address{Oakland University. mailing adress:\newline Mathematics and Science Center
Michigan.} 

\email{jorgedemellojr@oakland.edu}
\date{20 November 2026}

\begin{abstract}
Let $K$ be a number field and let $f\in\OK[x]$ be a monic polynomial of degree $d=p^r$ with $f\equiv x^d+c\pmod p$, such
as $x^2+c$. We prove that every $\alpha\in\Kbar$ of degree $D$ over $K$ with infinite orbit under $f$ satisfies
$$\hh_f(\alpha)\ \ge\ \frac{c_1}{D}\exp\bigl(-c_2\sqrt{\log3D}\bigr),$$
where $c_1,c_2>0$ depend only on $f$ and $K$. This implies $\hh_f(\alpha)\ge D^{-1-o(1)}$, making progress towards the
dynamical Lehmer conjecture. The polynomial $f$ need not be
postcritically finite. For such polynomials, for instance $x^2+1$, no lower bound of polynomial decay in $D$ was previously
known. For postcritically finite polynomials such as
$x^2-1$, it improves the recent quadratic bound of Habegger and Schmidt. The method is motivated by Dobrowolski's proof on
Lehmer's problem.
\end{abstract}

\maketitle
\tableofcontents

\section{Introduction}

Let $K$ be a number field, $f\in K[x]$ a polynomial of degree $d\ge2$, and denote by
$\hh_f(\alpha)=\lim_nd^{-n}h(f^{(n)}(\alpha))$ its canonical height, where $h$ is the absolute logarithmic Weil height. By
Northcott's theorem, $\hh_f(\alpha)=0$ exactly when $\alpha$ is preperiodic, that is, when its orbit
$\alpha,f(\alpha),f^{(2)}(\alpha),\dots$ is finite. The dynamical Lehmer conjecture \cite[Conj.~3.25]{Sil07} asserts that the
canonical height of the other points is bounded below in terms of their degree alone: there is $c(f,K)>0$ such that
\begin{equation}\label{eq:lehmer}
\hh_f(\alpha)\ \ge\ \frac{c(f,K)}{D},\qquad D=[K(\alpha):K],
\end{equation}
for every $\alpha\in\Kbar$ that is not preperiodic.

For $f=x^d$ this is Lehmer's problem, and the classical bound of Dobrowolski \cite{Dob},
$h(\alpha)\gg D^{-1}(\log\log3D/\log3D)^3$, is still the best known up to constants. His proof rests on the congruences
$A(x^p)\equiv A(x)^p\pmod p$, that is, on the fact that $x\mapsto x^p$ lifts the Frobenius at every prime $p$, used for all
$p\le\log D$; this is the starting point of our method. For an arbitrary $f$, Baker \cite[Thm.~1.14]{Bak} proved, using
the energy of dynamical Green functions, that \eqref{eq:lehmer} holds outside an exceptional set of size $O(D\log D)$.

For other maps, bounds of polynomial decay valid for all points were obtained only recently. In a key result, adapting Dimitrov's proof of
the Schinzel--Zassenhaus conjecture \cite{Dim}, Habegger and Schmidt proved the quadratic bound
\begin{equation}\label{eq:HS}
\hh_f(\alpha)\ \ge\ \frac c{D^2}
\end{equation}
for $x^p+c$ with $0$ periodic, under a technical hypothesis verified for instance when $p=2$ and $c$ is real \cite{HS21},
and more generally for centered, hyperbolic, postcritically finite polynomials of prime-power degree with integral
coefficients \cite{HS26}. Here a polynomial is \emph{postcritically finite} if all its critical points are preperiodic. The
integrality hypotheses of \cite{HS26}, in the spirit of Epstein's integrality theorem for postcritically finite polynomials
\cite{Eps}, imply that $f$ lifts a Frobenius at $p$ up to a translation \cite[Lemma~5.7]{HS26}; this is the property we use.
All these bounds concern postcritically finite polynomials, and, as far as we know, no bound of polynomial decay was known
for a polynomial with a critical point of infinite orbit, such as $x^2+1$. We show that a single prime at which $f$ lifts
Frobenius suffices to obtain the exponent $1+o(1)$ in \eqref{eq:lehmer} for all points, as follows.

\begin{definition}\label{def:shape}
Let $p$ be a prime and $d=p^r$. A monic $f\in\OK[x]$ of degree $d$ has \emph{Frobenius shape at $p$} if
$$f\ \equiv\ (x+b)^d+c\pmod{p\OK[x]}\qquad\text{for some }b,c\in\OK.$$
\end{definition}

Equivalently $f\equiv x^d+c'\pmod{p\OK[x]}$, with $c'=b^d+c$. Every $x^2+c$ with $c\in\OK$ has Frobenius shape at $2$, and so has every
monic $f\in\ZZ[x]$ of prime degree $p$ whose coefficients of $x,\dots,x^{p-1}$ are divisible by $p$, as well as the
polynomials of Habegger and Schmidt.

A point $x$ is \emph{periodic} if $f^{(n)}(x)=x$ for some $n\ge1$, and the finite critical points of $f$ are the zeros of $f'$.
Fix a prime $\p\mid p$ of $K$. For $\alpha\in\Kbar$ let $e_\p(\alpha)$ be the largest absolute ramification index of a prime
of $K(\alpha)$ above $\p$, so that $1\le e_\p(\alpha)\le e(\p|p)\,[K(\alpha):K]$, and write $\alpha_n=f^{(n)}(\alpha)$ and
$D_n=[K(\alpha_n):K]$.

\begin{mainthm}\label{thm:main}
Let $f\in\OK[x]$ be a monic polynomial of Frobenius shape at $p$ none of whose finite critical points is
periodic. There are $c_1,c_2>0$ such that every $\alpha\in\Kbar$ that is not preperiodic, with
$D=[K(\alpha):K]\ge2$, satisfies
\begin{enumerate}[label=\rm(\roman*),leftmargin=*]
\item $\displaystyle\hh_f(\alpha)\ \ge\ \frac{c_1}{e_\p(\alpha)\,D}\exp\bigl(-c_2\sqrt{\log3D}\bigr)$;
\item $\displaystyle\hh_f(\alpha)\ \ge\ \frac{c_1}{D}\exp\bigl(-(c_2+1)\sqrt{\log3D}\bigr)$.
\end{enumerate}
In particular $\hh_f(\alpha)\ge D^{-1-o(1)}$ for all $\alpha$.
\end{mainthm}

The constants depend only on $K$, $f$ and $\p$, and the polynomials covered include $x^2+c$ for every $c\in\OK$ such that $0$ is
not periodic (over $\ZZ$: every $c\ne0,-1$), and $x^d+c$ under the same condition (\S\ref{sec:examples}). The exponent $1$ in part (ii) is
optimal: if $\beta\in K$ is not preperiodic and $\alpha\in f^{-N}(\beta)$ has degree $d^N$, then $\hh_f(\alpha)=d^{-N}\hh_f(\beta)\asymp D^{-1}$.
Part (i) is the better bound when $e_\p(\alpha)<\exp(\sqrt{\log3D})$.

In \S\S\ref{sec:prelim}--\ref{sec:proofA}, among other things, we adapt Dobrowolski's method to a single prime. The characteristic polynomials of
the $\alpha_n$ satisfy a $p$-adic congruence whose modulus grows with $n$, and we compare it with an archimedean estimate in
a basis adapted to $f$ (the polynomials $\prod_j(f^{(j)})^{\varepsilon_j}$, with $\varepsilon_j$ the base-$d$ digits of the degree) over
$\asymp\sqrt{\log D}$ consecutive steps. Write $k_n=D_n/D_{n+1}\in\{1,\dots,d\}$; the degree \emph{drops} at step $n$ if $k_n\ge2$.
If the degree does not drop along $\asymp\sqrt{\log D}$ consecutive steps, we obtain the following bound for iterates of constant degree.

\begin{thmnn}
There are $c_1,c_2>0$ and $k(D)=O(\sqrt{\log3D})$ such that, for every non-preperiodic $\alpha\in\Kbar$ with $D\ge2$: if
$D_n$ is constant for $j\le n\le j+k(D)$, then
$$\hh_f(\alpha)\ \ge\ \frac{c_1}D\exp\bigl(-c_2\sqrt{\log3D}\bigr)\prod_{n<j}\frac{k_n}d.$$
\end{thmnn}

For integral $\alpha$ this is Proposition \ref{thm:window}; for the others it follows from Lemma \ref{lem:nonintegral}, which
gives $\hh_f(\alpha)\ge\log2/([K:\QQ]D)$.

Every step $n<j$ with $k_n<d$ lowers this bound by the factor $k_n/d$, so it is necessary to locate the steps at which the degree drops. Fix $\Kbar\subset\CC_p$ and the valuation $v$ with $v(p)=1$. A \emph{cycle} of $f$ is the orbit of a periodic point, and it is
\emph{superattracting} if it contains a critical point; let $\Sfin\subset\overline{\FF}_p$ be the set of reductions of the
points of the superattracting cycles (empty if no finite critical point is periodic), and let $\nu_\p\ge1$ be the integer part of the largest
valuation $v(f'(\eta))$ over the periodic points $\eta$ with $\bar\eta\notin\Sfin$ (Definition \ref{def:invariants}).

\begin{mainthm}\label{thm:droploc}
Let $f\in\OK[x]$ be a monic polynomial of Frobenius shape at $p$, and let $\alpha\in\Kbar$ be not preperiodic, with a
$K$-conjugate $\alpha'$ such that $v(\alpha')\ge0$ and $\bar\alpha'\notin\Sfin$. For every field $K\subseteq F\subseteq\Kbar$, if
$[F(\alpha_{n+1}):F]<[F(\alpha_n):F]$ then $n\le\lceil\log_de_\p(\alpha)\rceil+\nu_\p-1$.
\end{mainthm}

If the degree drops at step $n$, two conjugates of $\alpha_n$ have the same image under $f$. Both lie in the residue disc (the
set of integral points with a given reduction modulo the maximal ideal) of some periodic point $\eta$ and are $\p$-adically close to
$\eta$, and a divided difference shows that $v(f'(\eta))$ is large; this bounds $n$ (\S\ref{sec:drops}). Choosing $j$ after the last
drop, the product in the bound above is at least $d^{-j}\gg e_\p(\alpha)^{-1}$, which gives part (i) of Theorem \ref{thm:main}.
The same argument, applied to every fibre of $f$, shows that a drop at step $m$ forces at least half of the conjugates $\alpha'$ of $\alpha$
to have the property that $K_\p(\alpha')$ has ramification index $>d^{\,m-\nu_\p}$ over $\QQ_p$ (\S\ref{sec:clusters}). These
conjugates form $\p$-adic clusters, sets of $e/e(\p|p)$ pairwise $\p$-adically close conjugates, where $e$ is the ramification
index of $K_\p(\alpha')$, so the
discriminant of the minimal polynomial of $\alpha_m$ is divisible by a large power of $p$. This gives $\hh_f(\alpha)\gg D^{-1}$
unless $D_md^{\,m}\ll D\log D$, and in that case the bound above with $j=m+1$ gives $\hh_f(\alpha)\gg(D_md^{\,m})^{-1}$, up to the
factor $\exp(-c_2\sqrt{\log3D})$. Since $\log3D\le\exp(\sqrt{\log3D})$, the remaining factor $\log3D$ is absorbed into the exponential, and this proves part (ii). For $x^2+1$ the constants in part (i) can be made explicit (\S\ref{sec:explicit}), and so can those of Theorem \ref{thm:attract} for $x^2-1$ (Proposition \ref{prop:explicitminus}).

\begin{mainthm}\label{thm:B}
For $f=x^2+1$ and every $\alpha\in\Qbar$ that is not preperiodic, with $D=[\QQ(\alpha):\QQ]\ge2$ and $e_2(\alpha)$ the
ramification index $e_\p(\alpha)$ for $\p=2$,
$$\hh_f(\alpha)\ \ge\ \frac1{32\,e_2(\alpha)\,D}\exp\bigl(-8.4\sqrt{\log3D}\bigr).$$
\end{mainthm}

Theorem \ref{thm:droploc} does not apply to integral points all of whose conjugates reduce into $\Sfin$. For these, we use
instead that $f$ contracts $\p$-adically near a superattracting cycle. The valuation of the distance to the cycle is multiplied by $d$ at
each step until it reaches $1$, and after that it is multiplied by at least $2$ along each period of the cycle; a single resultant
then gives a bound that depends only on the ramification (\S\ref{sec:attract}).

\begin{definition}\label{def:omega}
Suppose $f$ has a superattracting cycle in $\Kbar$. For such a cycle $Z$ let $m_Z=\prod_c\deg_cf$, the product of the
local degrees $\deg_cf=1+\mathrm{ord}_{x=c}f'(x)$ of $f$ at the critical points $c\in Z$. Let $Q$ be the least common multiple of the periods of the
superattracting cycles, $M$ the total number of their points, $m_f=\min_Zm_Z\ge2$, and
$$\omega_f=Q\,\frac{\log d}{\log m_f}\ \ge\ 1,$$
the inequality because $m_f\le d^Q$. Finally let $C_f\ge0$ be such
that $|h(x)-\hh_f(x)|\le C_f$ for all $x\in\Kbar$ \cite[Thm.~3.20]{Sil07}.
\end{definition}
When there are superattracting cycles, we obtain the following
\begin{mainthm}\label{thm:attract}
Let $f$ be as in Theorem \ref{thm:droploc}, with a superattracting cycle in $\Kbar$. If $\alpha\in\Kbar$ is integral over
$\OK$, not preperiodic, and $\bar\alpha'\in\Sfin$ for every $K$-conjugate $\alpha'$ of $\alpha$, then
$$\hh_f(\alpha)\ \ge\ \frac B{2d^{2Q+1}}\min\Bigl\{1,\bigl(2AB\bigr)^{-\omega_f}\Bigr\}\cdot\frac1{e_\p(\alpha)},\qquad
A=\frac{M\,[K:\QQ]}{[K_\p:\QQ_p]\log p},\quad B=\log2+2C_f.$$
\end{mainthm}
Putting the results above together enables us to prove the Lehmer-type bound with exponent $1$ for all polynomials of Frobenius shape.
\begin{mainthm}\label{thm:all}
Let $f\in\OK[x]$ be any monic polynomial of Frobenius shape at $p$. There are $c_1,c_2>0$ such that every
non-preperiodic $\alpha\in\Kbar$ with $D\ge2$ satisfies
$$\hh_f(\alpha)\ \ge\ \frac{c_1}{D}\exp\bigl(-c_2\sqrt{\log3D}\bigr).$$
In particular $\hh_f(\alpha)\ge D^{-1-o(1)}$ for all such $\alpha$.
\end{mainthm}

Thus, for $x^2-1$, $x^2+c$ with $0$ periodic, and $x^4+2x^2$ at $p=2$, Theorem \ref{thm:all} gives $D^{-1-o(1)}$; for $x^2-1$ this
improves the bound $D^{-2}$ of \cite{HS21}. For $f=x^d$ it gives $h(\alpha)\ge D^{-1-o(1)}$, which is weaker than Dobrowolski's
bound. Here $f$ need not be postcritically finite.

In Sections \ref{sec:prelim}--\ref{sec:proofA} we prove the bound for iterates of constant degree (Proposition \ref{thm:window}). Section \ref{sec:local}
contains the local analysis at $\p$, and Sections \ref{sec:drops}, \ref{sec:clusters} and \ref{sec:attract} prove the results on
the location of the drops, the ramified clusters and the superattracting classes, respectively; Section \ref{sec:examples}
contains examples with explicit constants.


\section{Polynomials: setting and preliminaries}\label{sec:prelim}

In this section we fix the notation for polynomials and their canonical local heights, and prove two elementary estimates.

Throughout, $K$ is a number field, $p$ a prime, $d=p^r$, and $f\in\OK[x]$ is monic of
degree $d$. For an embedding $\sigma\colon K\to\CC$ write $\sigma f$ for the polynomial
with coefficients conjugated by $\sigma$, and for $F\in\CC[x]$ monic of degree $d$ let
$$\lambda_F(z)=\lim_{n\to\infty}d^{-n}\log\max\{1,|F^{(n)}(z)|\}.$$
Then $\lambda_F\ge0$ and $\lambda_F(F(z))=d\lambda_F(z)$. If $\alpha\in\Kbar$ is integral
over $\OK$, then all non-archimedean local heights of the orbit vanish, and
\begin{equation}\label{eq:hsum}
\hh_f(\alpha)=\frac1{[K(\alpha):\QQ]}\sum_{\tau}\lambda_{\tau|_Kf}(\tau\alpha),
\end{equation}
the sum running over the embeddings $\tau\colon K(\alpha)\to\CC$
\cite[\S5.9]{Sil07}, \cite{CG}.

We first bound the absolute value of a complex number in terms of its canonical local height.

\begin{lemma}\label{lem:size}
Let $F=x^d+\sum_{i<d}a_ix^i\in\CC[x]$ and $R=1+\sum_{i<d}|a_i|$. Then
$|w|\le R\,e^{\lambda_F(w)}$ for all $w\in\CC$. For $F=x^2+c$ one may take instead
$R=R_c=\frac12(1+\sqrt{1+4|c|})$.
\end{lemma}

\begin{proof}
If $|w|\le R$ this holds because $\lambda_F\ge0$. Let $|w|>R$ and $u=|w|/R>1$. Since
$|w|>R\ge1$, we have $\sum_{i<d}|a_i||w|^i\le(R-1)|w|^{d-1}$, so
$$|F(w)|\ \ge\ |w|^{d-1}\bigl(|w|-(R-1)\bigr).$$
Now $|w|-(R-1)\ge|w|R^{1-d}$: this is equivalent to $|w|(1-R^{1-d})\ge R-1$, which holds
because $|w|(1-R^{1-d})\ge R(1-R^{1-d})=R-R^{2-d}\ge R-1$, as $R\ge1$ and $d\ge2$. Hence
$|F(w)|\ge|w|^dR^{1-d}$, i.e.\ $|F(w)|/R\ge u^d>1$. By induction
$|F^{(n)}(w)|/R\ge u^{d^n}$; taking logarithms, dividing by $d^n$ and letting
$n\to\infty$ gives $\lambda_F(w)\ge\log u$.

For $F=x^2+c$ and $|w|>R_c$: $|F(w)|\ge|w|^2-|c|=R_cu^2+|c|(u^2-1)\ge R_cu^2$, using
$R_c^2=R_c+|c|$, and the same induction applies.
\end{proof}

We write $R_f=\max_\sigma R(\sigma f)$, with $R(\cdot)$ as in Lemma \ref{lem:size}.

Points that are not integral have canonical height at least of order $1/D$.

\begin{lemma}\label{lem:nonintegral}
If $\alpha\in\Kbar$ is not integral over $\OK$, then
$\hh_f(\alpha)\ge\log2/[K(\alpha):\QQ]$.
\end{lemma}

\begin{proof}
Normalise $|\cdot|_v$ to extend the usual absolute value on $\QQ$, so $|\ell|_v=\ell^{-1}$
for $v\mid\ell$, and for each place $v$ of $K(\alpha)$ let $n_v=e_vf_v$ and
$\lambda_v(\alpha)=\lim_nd^{-n}\log\max\{1,|f^{(n)}(\alpha)|_v\}$, so that $\hh_f(\alpha)=\sum_vn_v\lambda_v(\alpha)/[K(\alpha):\QQ]$. Let $v$ be a
place of $K(\alpha)$ above a prime $\ell$ with $|\alpha|_v>1$. Since $f$ is monic with
integral coefficients, $|f(\alpha)|_v=|\alpha|_v^d$, so
$\lambda_v(\alpha)=\log|\alpha|_v$; and $|\alpha|_v$ is an integral power of
$\ell^{1/e_v}$, so $\lambda_v(\alpha)\ge\frac1{e_v}\log\ell$. All local heights are
non-negative, so
$\hh_f(\alpha)\ge\frac{e_vf_v}{[K(\alpha):\QQ]}\lambda_v(\alpha)\ge\frac{\log2}{[K(\alpha):\QQ]}$.
\end{proof}

\section{The growing congruence}\label{sec:congruence}

In this section we show that the characteristic polynomials of the iterates of a point satisfy congruences modulo powers of $p$ whose exponent grows along the orbit.

For a monic $A\in\OK[X]$ with roots $y_1,\dots,y_D$ and $g\in\OK[x]$ put
$A^{[g]}(X)=\prod_i(X-g(y_i))\in\OK[X]$. Then $(A^{[g]})^{[h]}=A^{[h\circ g]}$.
Congruences between polynomials are coefficientwise, modulo ideals of $\OK$.

For $f=x^p$ the congruences below are the classical Euler--Fermat congruences for traces of powers
\cite{Smy86,Zar,Vin,BGW}; we need them for iterates of a polynomial of Frobenius shape.

Modulo $p$, the iterates of $f$ are powers of $x$ shifted by constants.

\begin{lemma}\label{lem:reduction}
Suppose $f$ has Frobenius shape at $p$, and put $c'=b^d+c$. Define $\gamma_0=0$ and
$\gamma_{n+1}=\gamma_n^d+c'$. Then $f^{(n)}(x)\equiv x^{d^n}+\gamma_n\pmod{p\OK[x]}$ for all
$n\ge0$.
\end{lemma}

\begin{proof}
Since $p$ divides $\binom di$ for $0<i<d$, we have $(x+b)^d\equiv x^d+b^d$, so
$f\equiv x^d+c'$. If $f^{(n)}\equiv x^{d^n}+\gamma_n$, then
$f^{(n+1)}\equiv(x^{d^n}+\gamma_n)^d+c'\equiv x^{d^{n+1}}+\gamma_n^d+c'$.
\end{proof}

Since $\OK/p\OK$ is finite, these constants and the corresponding power maps are eventually periodic.

\begin{lemma}\label{lem:stable}
There are integers $n_0\ge0$ and $q\ge1$, depending only on $K$, $p$ and $f$ modulo
$p$, such that for all $n\ge n_0$:
$$a^{d^{n+q}}\equiv a^{d^n}\pmod{p\OK}\quad\text{for all }a\in\OK,\qquad
\gamma_{n+q}\equiv\gamma_n\pmod{p\OK}.$$
\end{lemma}

\begin{proof}
Let $\mathcal R=\OK/p\OK$, a finite ring, and let $\phi_n\colon\mathcal R\to\mathcal R$ be
$a\mapsto a^{d^n}$. The pair $(\phi_n,\bar\gamma_n)$ lies in the finite set
$\mathcal R^{\mathcal R}\times\mathcal R$, and $(\phi_{n+1},\bar\gamma_{n+1})$ is
determined by $(\phi_n,\bar\gamma_n)$, since $\phi_{n+1}(a)=\phi_n(a)^d$ and
$\bar\gamma_{n+1}=\bar\gamma_n^d+\bar c'$. So the sequence is eventually periodic.
\end{proof}

We illustrate the integers $n_0$ and $q$ in two cases.

\begin{example}\label{ex:n0q}
\begin{enumerate}[label=(\alph*),leftmargin=*]
\item $K=\QQ$, $f=x^p+c$ with $c\in\ZZ$. On $\FF_p$, $a^p=a$, so $n_0=0$, and
$\gamma_n\equiv nc$; thus $q=p$, or $q=1$ if $p\mid c$. For $p=2$ and $c$ odd, $q=2$.
\item $K=\QQ(i)$, $f=x^2+i$, $p=2$. Here $\OK/2\OK=\FF_2[\varepsilon]/(\varepsilon^2)$ with
$\varepsilon=1+i$, and $(u+v\varepsilon)^2=u^2$; so $a^{2^n}\equiv a^{2^{n+1}}$ for $n\ge1$
but not for $n=0$. The orbit $\gamma_n$ is $0,i,i+1,i,i+1,\dots$ modulo $2$. Thus
$n_0=1$ and $q=2$.
\end{enumerate}
\end{example}

Modulo $p$, raising the roots of a polynomial to the $p$-th power raises its coefficients to the $p$-th power.

\begin{lemma}\label{lem:universal}
Let $e_j$ be the $j$-th elementary symmetric polynomial in $D$ variables. For each $j$ there is $S_j\in\ZZ[e_1,\dots,e_D]$ with
$e_j(y_1^p,\dots,y_D^p)=e_j(y)^p+pS_j(e(y))$.
\end{lemma}

\begin{proof}
Modulo $p$, the $p$-th power map is additive, so $e_j(y^p)-e_j(y)^p$ is a symmetric
polynomial with coefficients divisible by $p$. Divide by $p$ and apply the fundamental
theorem on symmetric polynomials.
\end{proof}

The following lemma is the main step: replacing the roots of two congruent polynomials by their images under $f$ raises the modulus of the congruence by one power of $p$.

\begin{lemma}\label{lem:raise}
Suppose $f$ has Frobenius shape at $p$. Let $B,B'\in\OK[X]$ be monic of the same degree
with $B\equiv B'\pmod{p^t\OK}$, $t\ge1$. Then $B^{[f]}\equiv B'^{[f]}\pmod{p^{t+1}\OK}$.
\end{lemma}

\begin{proof}
Write $f(x)=(x+b)^d+c+p\,g(x)$ with $g\in\OK[x]$, and let $y_i$, $y'_i$ be the roots of
$B$, $B'$, with elementary symmetric functions $e\equiv e'\pmod{p^t}$. Put $z_i=(y_i+b)^d+c$ and
$w_i=z_i+p\,g(y_i)=f(y_i)$, and similarly $z'_i$, $w'_i$; we first compare the polynomials with roots $z_i$ and $z'_i$.

Translation of roots by $b$ or $c$ changes the
coefficients by integral polynomials in them, so preserves congruences. Raising roots to
the $p$-th power raises the modulus: by Lemma \ref{lem:universal} the new coefficients
are $e_j^p+pS_j(e)$, and $pS_j(e)\equiv pS_j(e')\pmod{p^{t+1}}$, while
$e_j^p-e_j'^{\,p}=\sum_{i\ge1}\binom pie_j'^{\,p-i}(e_j-e'_j)^i$ is divisible by
$p^{t+1}$: for $0<i<p$ because $p\mid\binom pi$, and for $i=p$ because $tp\ge t+1$.
Doing this $r$ times, the polynomials with roots $z_i$ and $z'_i$ are congruent modulo
$p^{t+1}$ (even $p^{t+r}$).

Next we pass to the roots $w_i$. Expanding,
$$e_k(w)=e_k(z)+\sum_{j\ge1}p^j\sum_{\substack{S\subseteq\{1,\dots,D\}\\|S|=k}}
\ \sum_{\substack{T\subseteq S\\|T|=j}}\ \prod_{i\in T}g(y_i)\prod_{i\in S\setminus T}z_i.$$
For each $j$ the inner double sum is a symmetric polynomial in $y_1,\dots,y_D$ with
coefficients in $\OK$, since $z_i$ and $g(y_i)$ are polynomials in $y_i$ alone. By the
fundamental theorem on symmetric polynomials over the ring $\OK$, it equals
$\Psi_{k,j}(e(y))$ for a polynomial $\Psi_{k,j}\in\OK[e_1,\dots,e_D]$ that depends only
on $D$, $k$, $j$, $b$, $c$ and $g$, not on the roots. The same $\Psi_{k,j}$ therefore
serves for $B$ and $B'$, so $\Psi_{k,j}(e)\equiv\Psi_{k,j}(e')\pmod{p^t}$, and the terms
$p^j\Psi_{k,j}$ agree modulo $p^{t+1}$. Together with the first part, $e_k(w)\equiv e_k(w')\pmod{p^{t+1}}$.
\end{proof}

Iterating Lemma \ref{lem:raise}, we obtain the growing congruence.

\begin{proposition}\label{prop:congruence}
Suppose $f$ has Frobenius shape at $p$, and let $n_0,q$ be as in Lemma
\ref{lem:stable}. For every monic $A\in\OK[X]$, the polynomials $A_n=A^{[f^{(n)}]}$
satisfy
$$A_{n+q}\equiv A_n\pmod{p^{\,n-n_0+1}\OK}\qquad(n\ge n_0).$$
Consequently $A_m\equiv A_n\pmod{p^{\,n-n_0+1}}$ whenever $m\ge n\ge n_0$ and
$q\mid m-n$.
\end{proposition}

\begin{proof}
Let $y_i$ be the roots of $A$ and $\OO$ the ring of integers of a field containing them.
By Lemma \ref{lem:reduction}, $f^{(n)}(y_i)=y_i^{d^n}+\gamma_n+p\,\eta_i$ with
$\eta_i\in\OO$. Hence, modulo $p\OO$, the coefficients of $A_n$ are those of the
polynomial with roots $y_i^{d^n}+\gamma_n$. Applying Lemma \ref{lem:universal}
$rn$ times gives, by induction, $e_j(y^{d^n})\equiv e_j(y)^{d^n}\pmod{p\OO}$. The
coefficients of $\prod_i(X-u_i-\gamma)$ are given by one fixed $\ZZ$-polynomial in
$e(u)$ and $\gamma$, independent of $n$. So modulo $p\OO$ the coefficients of $A_n$ are
this fixed polynomial evaluated at $e(y)^{d^n}$ and $\gamma_n$. By Lemma \ref{lem:stable},
$e(y)^{d^{n_0+q}}\equiv e(y)^{d^{n_0}}$ and $\gamma_{n_0+q}\equiv\gamma_{n_0}$ modulo $p$,
since $e(y)\in\OK$. So $A_{n_0+q}\equiv A_{n_0}\pmod{p\OO}$, and as both sides lie in
$\OK[X]$ and $p\OO\cap\OK=p\OK$, the congruence holds modulo $p\OK$.

Applying Lemma \ref{lem:raise} $n-n_0$ times gives
$$A_{n+q}=A_{n_0+q}^{[f^{(n-n_0)}]}\equiv A_{n_0}^{[f^{(n-n_0)}]}=A_n\pmod{p^{\,n-n_0+1}}.$$
The last statement follows by chaining, as all the moduli
$p^{\,n-n_0+1},p^{\,n+q-n_0+1},\dots$ are multiples of $p^{\,n-n_0+1}$.
\end{proof}

For minimal polynomials, the growing congruence gives a lower bound for the $p$-adic valuation of resultants.

\begin{corollary}\label{cor:resultant}
Let $\alpha$ be integral over $\OK$ of degree $D$ over $K$, let $P_n$ be the minimal
polynomial of $f^{(n)}(\alpha)$ over $K$, and $D_n=\deg P_n$. If $n_0\le n<m$, $q\mid m-n$
and $P_n\ne P_m$, then
$$v_p\bigl(\Nm_{K/\QQ}\Res(P_n,P_m)\bigr)\ \ge\ [K:\QQ]\,(n-n_0+1)\,\frac{D_nD_m}D.$$
\end{corollary}

\begin{proof}
Let $A$ be the minimal polynomial of $\alpha$ over $K$, so $A_k=P_k^{D/D_k}$. Write
$A_m=A_n+p^{\,n-n_0+1}B$ with $B\in\OK[X]$. For each root $y$ of $P_n$, $A_m(y)=p^{\,n-n_0+1}B(y)$,
so $\Res(P_n,A_m)=p^{(n-n_0+1)D_n}\Res(P_n,B)$ with $\Res(P_n,B)\in\OK$. Taking norms,
$p^{[K:\QQ](n-n_0+1)D_n}$ divides $\Nm\Res(P_n,A_m)=\Nm\Res(P_n,P_m)^{D/D_m}$.
\end{proof}

\section{Degrees along the orbit, and the adapted basis}\label{sec:basis}

In this section we study how the degrees behave along the orbit, and introduce a basis of polynomials adapted to $f$, whose values are controlled by the canonical local heights.

\begin{lemma}\label{lem:degrees}
Let $\alpha\in\Kbar$ have infinite $f$-orbit, $D=[K(\alpha):K]$, and
$D_n=[K(f^{(n)}(\alpha)):K]$. Then:
\begin{enumerate}[label=\rm(\alph*)]
\item $D_{n+1}=D_n/k_n$ with $k_n=[K(f^{(n)}(\alpha)):K(f^{(n+1)}(\alpha))]\le d$, and
$D_n\mid D$. Hence $D/D_n$ divides $\delta_d(D)$, the largest divisor of $D$ all of whose prime factors are
$\le d$, and $D_n\ge D/\delta_d(D)$ for all
$n$.
\item For $n\ne m$, no $K$-conjugate of $f^{(n)}(\alpha)$ equals one of
$f^{(m)}(\alpha)$; in particular $P_n\ne P_m$.
\end{enumerate}
\end{lemma}

\begin{proof}
(a) Put $y=f^{(n)}(\alpha)$. Then $K(f(y))\subseteq K(y)\subseteq K(\alpha)$, and $y$ is a
root of $f(X)-f(y)$, of degree $d$, so $k_n\le d$. Thus $D/D_n=k_0\cdots k_{n-1}$ divides
$D$ and has all prime factors $\le d$, so it divides $\delta_d(D)$.

(b) $\hh_f$ is invariant under conjugation and $\hh_f(f^{(k)}(\alpha))=d^k\hh_f(\alpha)$,
with $\hh_f(\alpha)>0$.
\end{proof}

For $m\ge0$ with base-$d$ expansion $m=\sum_j\varepsilon_jd^j$, $0\le\varepsilon_j<d$, put
$$B_m=\prod_j\bigl(f^{(j)}\bigr)^{\varepsilon_j}\in\OK[x],\qquad\sigma_d(m)=\sum_j\varepsilon_j.$$

The polynomials $B_m$ are monic of degree $m$, and their values are controlled by the canonical local height.

\begin{lemma}\label{lem:basis}
$B_m$ is monic of degree $m$, and for every embedding $\sigma$ and $z\in\CC$,
$$|\sigma B_m(z)|\ \le\ R_f^{\,\sigma_d(m)}\,e^{m\lambda_{\sigma f}(z)}.$$
Moreover $\sigma_d(m)\le(d-1)(\log_dN+1)$ for $0\le m<N$.
\end{lemma}

\begin{proof}
$B_m$ is a product of monic polynomials of total degree $\sum\varepsilon_jd^j=m$. By Lemma
\ref{lem:size} for $w=\sigma f^{(j)}(z)$,
$|\sigma f^{(j)}(z)|\le R_fe^{d^j\lambda_{\sigma f}(z)}$; multiply. An integer $m<N$
has at most $\log_dN+1$ digits, each at most $d-1$.
\end{proof}

\section{A bound for iterates of constant degree}\label{sec:proofA}

In this section we prove the bound for iterates of constant degree by comparing a $p$-adic lower bound and an archimedean upper bound for a Vandermonde determinant built from several consecutive points of the orbit.

Throughout, $n_0$ and $q$ are as in Lemma
\ref{lem:stable}, $R_f$ is as defined after Lemma \ref{lem:size}, and for $\beta\in\Kbar$ we write $\beta_n=f^{(n)}(\beta)$.

\subsection{Outline of the argument}
Let $\beta$ be an algebraic integer of degree $D'$ over $K$, and suppose that the degree of $\beta_n$ stays equal to $D'$ for
$n_0\le n\le n_0+L$. Collect all the conjugates of $\beta_{n_0},\beta_{n_0+1},\dots,\beta_{n_0+L}$: these are $N=(L+1)D'$
distinct algebraic integers, and the square of their Vandermonde determinant, multiplied over the embeddings of $K$, is
a non-zero rational integer. We estimate this integer in two ways.
\begin{enumerate}[label=(\roman*),leftmargin=*]
\item By the growing congruence (Corollary \ref{cor:resultant}), the conjugates of $\beta_n$ and
of $\beta_m$ are $p$-adically close whenever $q\mid m-n$, and closer the larger $n$ is: the resultant of their minimal
polynomials is divisible by $p^{(n-n_0+1)D'}$. There are about $L^2/(2q)$ such pairs of levels, with
weights $n-n_0+1$ up to $L$, and adding up these weights shows that, on average over the embeddings, the determinant has absolute
value at least about $p^{D'L^3/(6q)}$, since the integer is non-zero and divisible by a high power of $p$.
\item By Hadamard's inequality the determinant is at most the product of the lengths
of its rows. If the rows are written in the monomial basis $1,x,\dots,x^{N-1}$, the row of a point $z$ has length about
$\max(1,|z|)^N$, and the product of these over all points is controlled only by the Weil height, which may be much larger than the canonical height. Instead we use
the adapted basis $B_m$ of \S\ref{sec:basis}, which is triangular with respect to the monomials (so the determinant does
not change) and whose values are controlled by the canonical local height $\lambda$ (Lemma \ref{lem:basis}). Then the row
of $z$ has length about $e^{N\lambda(z)}$ times a factor $N^{O(1)}$, and the sum of $\lambda$ over all the points is
$D'\sum_nd^n\hh_f(\beta)\approx D'd^{\,L}\hh_f(\beta)$ (in this heuristic we ignore the fixed shift $n_0$). So the determinant is at most about
$\exp\bigl(O(N\log N)+N\cdot D'd^{\,L}\hh_f(\beta)\bigr)$.
\end{enumerate}
Comparing, $D'L^3\log p/(6q)\lesssim D'L\log N+(L+1)D'^2d^{\,L}\hh_f(\beta)$, that is,
$$\hh_f(\beta)\ \gtrsim\ \frac{L^2\log p/(6q)-\log N}{D'\,d^{\,L}}.$$
The numerator is positive as soon as $L^2\gg q\log N\approx q\log D'$, so the best choice is $L\asymp\sqrt{q\log D'}$,
which gives $\hh_f(\beta)\gg D'^{-1}d^{-O(\sqrt{q\log D'})}$. The rest of the section makes this precise.

\subsection{The determinant}
Fix $\beta$ integral over $\OK$ with infinite $f$-orbit, put $D'=[K(\beta):K]$, and let $L\ge1$ be an integer such that
\begin{equation}\label{eq:nodrop5}
[K(\beta_n):K]=D'\qquad\text{for }0\le n\le n_0+L.
\end{equation}
Then $K(\beta_n)=K(\beta)$ for these $n$, since $K(\beta_n)\subseteq K(\beta)$ and the degrees agree. Let $P_n$ be the minimal
polynomial of $\beta_n$ over $K$, monic with coefficients in $\OK$, and put $N=(L+1)D'$.

For an embedding $\sigma\colon K\to\CC$, let $Z_\sigma$ be the set of the numbers $\tau(\beta_n)$, where $n_0\le n\le n_0+L$ and
$\tau$ runs over the $D'$ embeddings $K(\beta)\to\CC$ extending $\sigma$. For fixed $n$ these are the $D'$ roots of $\sigma P_n$,
which are distinct, and by Lemma \ref{lem:degrees}(b) no root of $\sigma P_n$ is a root of $\sigma P_m$ for $n\ne m$. So
$|Z_\sigma|=N$. Order $Z_\sigma$ in any way, and put
$$\Delta_\sigma=\det\bigl(\sigma B_m(z)\bigr)_{z\in Z_\sigma,\ 0\le m<N},$$
where $\sigma B_m$ is the polynomial $B_m$ of \S\ref{sec:basis} with coefficients conjugated by $\sigma$.

The determinant $\Delta_\sigma$ is a Vandermonde determinant, so the product of the $\Delta_\sigma^2$ is a non-zero integer built from discriminants and resultants.

\begin{lemma}\label{lem:detint}
Up to sign, $\Delta_\sigma$ is the Vandermonde determinant $\prod_{z<z'}(z'-z)$ of $Z_\sigma$, and
$$\prod_\sigma\Delta_\sigma^2\ =\ \pm\,\Nm_{K/\QQ}\Bigl(\prod_{n_0\le n\le n_0+L}\disc(P_n)\prod_{n_0\le n<m\le n_0+L}\Res(P_n,P_m)^2\Bigr).$$ This is a non-zero rational integer, and each norm $\Nm\Res(P_n,P_m)$ is a non-zero
rational integer.
\end{lemma}

\begin{proof}
Since $\sigma B_m$ is monic of degree $m$ (Lemma \ref{lem:basis}), the matrix $(\sigma B_m(z))$ is obtained from the
Vandermonde matrix $(z^m)$ by adding to each column $m$ a linear combination of the columns $0,\dots,m-1$; these
operations do not change the determinant. Group the factors $z'-z$ of the Vandermonde determinant according to the levels of $z$ and $z'$. The
factors with both points in the level $n$ form a product whose square is $\pm\disc(\sigma P_n)$; the factors with one point in
the level $n$ and the other in the level $m>n$ form, up to sign, the product of $y'-y$ over the roots $y$ of $\sigma P_n$ and $y'$ of
$\sigma P_m$, which is $\pm\Res(\sigma P_n,\sigma P_m)$ because both polynomials are monic. Hence $\Delta_\sigma^2=\pm\sigma\bigl(\prod\disc(P_n)\prod\Res(P_n,P_m)^2\bigr)$,
and the product over $\sigma$ is the norm. The elements $\disc(P_n)$ and $\Res(P_n,P_m)$ lie in $\OK$ because the $P_n$ are
monic with integral coefficients, and they are non-zero because the points of $Z_\sigma$ are distinct.
\end{proof}

\subsection{The lower bound from the prime $p$}
Put
$$S_q(L)=\sum_{j=0}^L(j+1)\Bigl\lfloor\frac{L-j}q\Bigr\rfloor.$$
This counts, with weight $j+1$, the pairs of levels $n=n_0+j<m\le n_0+L$ with $q\mid m-n$: for fixed $j$ there are
$\lfloor(L-j)/q\rfloor$ such $m$.

By the growing congruence, this integer is divisible by a large power of $p$.

\begin{lemma}\label{lem:detlow}
Under \eqref{eq:nodrop5},
\begin{equation}\label{eq:Sq}
\sum_\sigma\log|\Delta_\sigma|\ \ge\ [K:\QQ]\,D'\,S_q(L)\log p,\qquad
S_q(L)\ \ge\ \frac{(L+1)(L+2)(L-3q+3)}{6q}.
\end{equation}
\end{lemma}

\begin{proof}
Let $\mathcal N=\prod_\sigma\Delta_\sigma^2$, a non-zero integer by Lemma \ref{lem:detint}. Then $|\mathcal N|\ge p^{v_p(\mathcal N)}$, that is,
$2\sum_\sigma\log|\Delta_\sigma|\ge v_p(\mathcal N)\log p$. All the factors in Lemma \ref{lem:detint} are integers, so $v_p(\mathcal N)$ is at
least the contribution of the resultants:
$$v_p(\mathcal N)\ \ge\ 2\sum_{\substack{n_0\le n<m\le n_0+L\\ q\mid m-n}}v_p\bigl(\Nm\Res(P_n,P_m)\bigr).$$
For each such pair, Corollary \ref{cor:resultant}, applied with $\alpha=\beta$ and $D_n=D_m=D=D'$, gives
$v_p(\Nm\Res(P_n,P_m))\ge[K:\QQ](n-n_0+1)D'$; its hypothesis $P_n\ne P_m$ holds by Lemma \ref{lem:degrees}(b). Writing
$n=n_0+j$ and summing over $m$ gives the first inequality.

For the second, note that for every integer $u\ge0$ one has $u-q\lfloor u/q\rfloor\le q-1$, so
$\lfloor u/q\rfloor\ge(u-q+1)/q$. Apply this with $u=L-j$, and use $\sum_{j=0}^L(j+1)(L-j)=\frac{L(L+1)(L+2)}6$ and
$\sum_{j=0}^L(j+1)=\frac{(L+1)(L+2)}2$:
\begin{align*}S_q(L)\ &\ge\ \frac1q\sum_{j=0}^L(j+1)(L-j)-\frac{q-1}q\sum_{j=0}^L(j+1)\\
&=\frac1q\cdot\frac{L(L+1)(L+2)}6-\frac{q-1}q\cdot\frac{(L+1)(L+2)}2,\end{align*}
which equals $\frac{(L+1)(L+2)(L-3q+3)}{6q}$.
\end{proof}

\subsection{The upper bound from the archimedean places}
Put
$$\Lambda=\Lambda(N)=\tfrac12\log N+(d-1)\Bigl(\frac{\log N}{\log d}+1\Bigr)\log R_f.$$

On the other hand, Hadamard's inequality and Lemma \ref{lem:basis} bound $\Delta_\sigma$ from above in terms of the canonical height of $\beta$.

\begin{lemma}\label{lem:detup}
Under \eqref{eq:nodrop5},
$$\sum_\sigma\log|\Delta_\sigma|\ \le\ [K:\QQ]\,N\Lambda+2[K:\QQ]\,N\,D'\,d^{\,n_0+L}\,\hh_f(\beta).$$
\end{lemma}

\begin{proof}
Fix $\sigma$. By Hadamard's inequality, $|\Delta_\sigma|$ is at most the product over $z\in Z_\sigma$ of the Euclidean lengths of the
rows $(\sigma B_m(z))_{0\le m<N}$. By Lemma \ref{lem:basis}, each entry satisfies
$$|\sigma B_m(z)|\ \le\ R_f^{\,\sigma_d(m)}e^{m\lambda_{\sigma f}(z)}\ \le\ R_f^{\,(d-1)(\log_dN+1)}\,e^{N\lambda_{\sigma f}(z)},$$
using $m<N$, $\lambda_{\sigma f}\ge0$ and $R_f\ge1$. A row has $N$ entries, so its length is at most $\sqrt N$ times this bound.
Taking logarithms and summing over the $N$ rows,
$$\log|\Delta_\sigma|\ \le\ N\Lambda+N\sum_{z\in Z_\sigma}\lambda_{\sigma f}(z).$$
Now sum over $\sigma$. For fixed $n$, the points $\tau(\beta_n)$, with $\sigma$ and $\tau$ varying, are the images of $\beta_n$ under
all the embeddings of $K(\beta_n)=K(\beta)$ into $\CC$; by \eqref{eq:hsum},
$$\sum_\sigma\sum_\tau\lambda_{\sigma f}(\tau\beta_n)=[K(\beta):\QQ]\,\hh_f(\beta_n)=[K:\QQ]\,D'\,d^{\,n}\hh_f(\beta),$$
because $\beta_n$ is integral and $\hh_f(f^{(n)}(\beta))=d^n\hh_f(\beta)$. Summing over $n_0\le n\le n_0+L$ and using
$\sum_{n\le n_0+L}d^n\le\frac d{d-1}d^{\,n_0+L}\le2d^{\,n_0+L}$ gives the claim.
\end{proof}

\subsection{Comparison and choice of the number of steps}

Comparing the two bounds and choosing $L$ of order $\sqrt{q\log D'}$, we obtain the bound when the degree is constant from the first step on, for a number of steps of order $\sqrt{q\log D'}$.

\begin{proposition}\label{prop:window}
There are $a\ge1$ and $c_3>0$, depending only on $K$ and $f$, with the following
property. Let $\beta\in\Kbar$ be integral over $\OK$ with infinite $f$-orbit, let
$D'=[K(\beta):K]$, let $\ell\ge\max\{\log6,\log3D'\}$, and put
$L=\lceil a\sqrt{q\ell}\,\rceil$. If $[K(f^{(n)}(\beta)):K]=D'$ for $0\le n\le n_0+L$, then
$$\hh_f(\beta)\ \ge\ \frac{c_3}{D'}\,d^{-a\sqrt{q\ell}}.$$
\end{proposition}

\begin{proof}
Combine Lemmas \ref{lem:detlow} and \ref{lem:detup}, divide by $[K:\QQ]D'$, and use $N=(L+1)D'$:
$$S_q(L)\log p\ \le\ (L+1)\Lambda+2(L+1)\,D'\,d^{\,n_0+L}\,\hh_f(\beta).$$
Equivalently,
\begin{equation}\label{eq:keyA}
\hh_f(\beta)\ \ge\ \frac{S_q(L)\log p-(L+1)\Lambda}{2(L+1)\,D'\,d^{\,n_0+L}}.
\end{equation}
This is useful only if the numerator is positive, and we choose $L$ so that it is at least $(L+1)\Lambda$, that is,
\begin{equation}\label{eq:claimA}
S_q(L)\log p\ \ge\ 2(L+1)\Lambda.
\end{equation}

The quantity $\Lambda$ grows like $\log N$: indeed $\Lambda\le b_1\log N+b_2$ with
$b_1=\frac12+(d-1)\log R_f/\log d$ and $b_2=(d-1)\log R_f$, which depend only on $d$ and $R_f$. By \eqref{eq:Sq}, the inequality \eqref{eq:claimA} follows
from
$$(L+2)(L-3q+3)\ \ge\ \frac{12q}{\log p}\bigl(b_1\log N+b_2\bigr).$$
We now bound both sides in terms of $\ell$. On the left, if $a\ge6\sqrt q/\sqrt{\log6}$, then $L\ge a\sqrt{q\log6}\ge6q$, so
$L-3q+3\ge L/2$ and $(L+2)(L-3q+3)\ge\frac12L^2\ge\frac12a^2q\ell$. On the right, $\log D'\le\ell$ and
$L+1\le a\sqrt q\,\ell+2$ (as $\sqrt\ell\le\ell$ for $\ell\ge1$), so
$$\log N=\log(L+1)+\log D'\ \le\ \ell+\log(a\sqrt q\,\ell+2)\ \le\ 2\ell+\log(a\sqrt q+2),$$
using $\log(a\sqrt q\,\ell+2)\le\log\ell+\log(a\sqrt q+2)\le\ell+\log(a\sqrt q+2)$. Hence \eqref{eq:claimA} holds as soon as
\begin{equation}\label{eq:achoice}
\tfrac12a^2\,\ell\ \ge\ \frac{12}{\log p}\Bigl(2b_1\,\ell+b_1\log(a\sqrt q+2)+b_2\Bigr)\qquad\text{for all }\ell\ge\log6.
\end{equation}
Dividing by $\ell\ge\log6$, it suffices that $\frac12a^2\ge\frac{12}{\log p}\bigl(2b_1+(b_1\log(a\sqrt q+2)+b_2)/\log6\bigr)$. The left
side grows like $a^2$ and the right side like $\log a$, so this holds for all large $a$, and we fix such an
$a\ge6\sqrt q/\sqrt{\log6}$, depending only on $K$ and $f$.

With this $a$, \eqref{eq:claimA} holds, so the numerator of \eqref{eq:keyA} is at least $(L+1)\Lambda$, and
$$\hh_f(\beta)\ \ge\ \frac{\Lambda}{2D'\,d^{\,n_0+L}}.$$
Since $N\ge2$, $\Lambda\ge\Lambda_0:=\frac12\log2$. Finally $L\le a\sqrt{q\ell}+1$, so $d^{\,n_0+L}\le d^{\,n_0+1}d^{\,a\sqrt{q\ell}}$, and the
proposition holds with $c_3=\Lambda_0/(2d^{\,n_0+1})$.
\end{proof}

\subsection{Iterates of constant degree from an arbitrary step}

Applying Proposition \ref{prop:window} to $f^{(j)}(\alpha)$, we obtain the bound when the degree is constant from an arbitrary step $j$ on.

\begin{proposition}\label{thm:window}
Let $f\in\OK[x]$ have Frobenius shape at $p$, let $\alpha\in\Kbar$ be integral over $\OK$ with infinite orbit,
$D=[K(\alpha):K]\ge2$, $D_n=[K(f^{(n)}(\alpha)):K]$ and $k_n=D_n/D_{n+1}\in\{1,\dots,d\}$. Let
$\ell=\log3D$, let $a,c_3$ be as in Proposition \ref{prop:window} and $L=\lceil a\sqrt{q\ell}\,\rceil$. If
$j\ge0$ is such that $D_n$ is constant for $j\le n\le j+n_0+L$, then
$$\hh_f(\alpha)\ \ge\ \frac{c_3}D\,d^{-a\sqrt{q\ell}}\prod_{n<j}\frac{k_n}d.$$
\end{proposition}

\begin{proof}
Put $\beta=f^{(j)}(\alpha)$ and $D'=D_j\le D$. Then $\beta$ is integral, as $f$ is monic with integral coefficients, and it has
infinite orbit. Since $D'\le D$ we have $\ell=\log3D\ge\max\{\log6,\log3D'\}$, and the degrees $[K(f^{(n)}(\beta)):K]=D_{j+n}$ are equal
to $D'$ for $0\le n\le n_0+L$. So Proposition \ref{prop:window} applies to $\beta$ with this $\ell$ and gives
$\hh_f(\beta)\ge c_3D'^{-1}d^{-a\sqrt{q\ell}}$. Now $\hh_f(\beta)=d^j\hh_f(\alpha)$, and $D/D'=k_0k_1\cdots k_{j-1}$ by Lemma
\ref{lem:degrees}(a). Hence
$$\hh_f(\alpha)\ =\ d^{-j}\hh_f(\beta)\ \ge\ \frac{c_3}{D'd^j}\,d^{-a\sqrt{q\ell}}\ =\ \frac{c_3}D\,d^{-a\sqrt{q\ell}}\prod_{n<j}\frac{k_n}d.$$
\end{proof}


\section{Local dynamics at \texorpdfstring{$\p$}{p}}\label{sec:local}

In this section we study the dynamics of $f$ on the residue discs of $\CC_p$: each disc contains exactly one periodic point, and orbits approach it at a rate that depends on the ramification.

Throughout, $f\in\OK[x]$ is a monic polynomial of degree $d=p^r$ with Frobenius shape at $p$ (Definition
\ref{def:shape}), and $\p\mid p$ is a prime of $K$. Let $\CC_p$ be the completion of an algebraic closure of $K_\p$,
fix an embedding $\Kbar\hookrightarrow\CC_p$ over $K$, let $v$ be the valuation with $v(p)=1$, $\OO_{\CC_p}$ the
valuation ring, and $x\mapsto\bar x$ the reduction $\OO_{\CC_p}\to\overline{\FF}_p$. For a polynomial $H$ over $\CC_p$ let
$v(H)$ be the least valuation of its coefficients.

We first write $f$ in a form adapted to its reduction modulo $p$.

\begin{definition}\label{def:normal}
Write $c'=b^d+c$, with $b,c$ as in Definition \ref{def:shape}, and let $h\in\OK[x]$ be the polynomial with
$$f(x)=x^d+c'+p\,h(x).$$
\end{definition}

Such an $h$ exists because $(x+b)^d\equiv x^d+b^d\pmod{p\OK[x]}$. Put
$$\rho(s)=\min\{ds,\ s+1\}\qquad(s\ge0),$$
an increasing function with $\rho(s)>s$ for $s>0$.

\subsection*{Residue discs}
For $x,y\in\OO_{\CC_p}$, $v(x-y)>0$ if and only if $\bar x=\bar y$. For $\xi\in\overline{\FF}_p$ let
$U_\xi=\{x\in\OO_{\CC_p}:\bar x=\xi\}$ be its residue disc, which is open and closed in $\CC_p$.

The first lemma describes $f$ on a residue disc, and introduces the invariant $\vartheta(x)=v(f'(x))$.

\begin{lemma}\label{lem:loc}
\begin{enumerate}[label=\rm(\alph*)]
\item The reduction $\bar f(x)=x^d+\bar c'$ is a bijection of $\overline{\FF}_p$, and every point of
$\overline{\FF}_p$ is periodic for it; let $q_\xi$ be the exact period of $\xi$.
\item For $x\in\OO_{\CC_p}$ let $a(t)=f(x+t)-f(x)$. Then $a(0)=0$ and
$$a(t)=t^d+t\,a_1(t),\qquad a_1\in\OO_{\CC_p}[t],\quad v(a_1)\ge1.$$
\item For $x,y\in\OO_{\CC_p}$ with $\bar x=\bar y$: $v(f(x)-f(y))\ge\rho(v(x-y))$.
\item In the notation of (b), put $\vartheta(x)=v(a'(0))=v(f'(x))$. One has $\vartheta(x)\ge1$,
$\vartheta(x)=\infty$ exactly when $x$ is a critical point of $f$, and there is $\vartheta_0$ with
$\vartheta(x)=\vartheta_0$ for all $x$ outside finitely many residue discs.
\item If $x\in\CC_p$ is not integral, then $v(f(x))=d\,v(x)<v(x)<0$.
\end{enumerate}
\end{lemma}

\begin{proof}
(a) Let $\mathbf k$ be the residue field of $\p$. The reduction is the composite of $x\mapsto x^{p^r}$ and the
translation by $\bar c'\in\mathbf k$, two bijections of $\overline{\FF}_p$ which preserve every finite field $\FF$
containing $\mathbf k$; on the finite set $\FF$ a bijection is a permutation, and every point of $\overline{\FF}_p$ lies
in such an $\FF$.

(b) By Definition \ref{def:normal}, $a(t)-t^d=\sum_{0<i<d}\binom di x^{d-i}t^i+p\bigl(h(x+t)-h(x)\bigr)$. Since
$p\mid\binom di$ for $0<i<d$, all coefficients of this polynomial in $t$ have valuation $\ge1$, and it vanishes at
$t=0$; so it equals $t\,a_1(t)$ as claimed.

(c) Put $t=y-x$, so that $v(t)=v(x-y)>0$ and $f(y)=f(x)+a(t)$. By (b),
$v(f(x)-f(y))=v(t^d+t\,a_1(t))\ge\min\{dv(t),\,v(t)+1\}$.

(d) Clearly $a'(0)=f'(x)$. By Definition \ref{def:normal}, $f'(x)=dx^{d-1}+p\,h'(x)$, and both terms have valuation
$\ge1$ for integral $x$. The polynomial $f'\in\OK[x]$ is non-zero, with leading coefficient $d$; write $f'=\lambda g_0$
with $\lambda\in K_\p^\times$ and $g_0$ having a unit coefficient. The reduction $\bar g_0\ne0$ has at most $d-1$ zeros
in $\overline{\FF}_p$, and outside their residue discs $\vartheta=v(\lambda)$.

(e) If $v(x)<0$, then every term $a_ix^i$, $i<d$, of $f$ has valuation $\ge iv(x)>dv(x)=v(x^d)$.
\end{proof}

Using the contraction in Lemma \ref{lem:loc}(c), we find exactly one periodic point in each residue disc.

\begin{lemma}\label{lem:periodic}
Each residue disc $U_\xi$ contains exactly one periodic point $\eta_\xi$ of $f$. It has exact period $q_\xi$,
$f(\eta_\xi)=\eta_{\bar f(\xi)}$, and if $x\in U_\xi$ is not periodic then $v(x-\eta_\xi)=v(f^{(q_\xi)}(x)-x)$. If
$F\subseteq\CC_p$ is a closed subfield containing $K_\p$ and $U_\xi\cap F\ne\emptyset$, then $\eta_\xi\in F$. Every
periodic point of $f$ in $\CC_p$ is integral, so the $\eta_\xi$, $\xi\in\overline{\FF}_p$, are all the periodic points.
\end{lemma}

\begin{proof}
Put $q=q_\xi$ and $G=f^{(q)}$. Then $G(U_\xi)\subseteq U_\xi$, and by Lemma \ref{lem:loc}(c),
$v(G(x)-G(y))>v(x-y)$ for $x\ne y$ in $U_\xi$; so $G$ has at most one fixed point in $U_\xi$. For
existence, take $x\in U_\xi$ (in $F$ if possible) and put $s_k=v(G^{k+1}(x)-G^k(x))$. Then $s_0>0$ and
$s_{k+1}\ge\rho(s_k)\ge s_k+\min\{(d-1)s_0,1\}$, so $s_k\to\infty$. By the ultrametric inequality the sequence
$G^k(x)$ is Cauchy; its limit $\eta$ lies in the closed set $U_\xi$, and in $F$ if $x$ does, as $f$ has coefficients
in $K$. It is fixed by $G$, which is continuous ($v(G(x)-G(y))\ge v(x-y)$). For $0<j<q$, $f^{(j)}(\eta)$
lies in $U_{\bar f^j(\xi)}\ne U_\xi$, so $\eta$ has exact period $q$. A periodic point in $U_\xi$
has period divisible by $q$, since its period fixes $\xi$ under $\bar f$; so it is a fixed point
of some $G^k$, $k\ge1$, and the argument above, applied to $G^k$, shows that it equals $\eta$. The
point $f(\eta_\xi)$ is periodic and lies in $U_{\bar f(\xi)}$, so it equals
$\eta_{\bar f(\xi)}$. If $x\in U_\xi$ is not periodic and $s=v(x-\eta_\xi)$, then
$v(G(x)-\eta_\xi)=v(G(x)-G(\eta_\xi))>s$, so $v(G(x)-x)=s$.

A point of period $q$ is a root of the monic polynomial $f^{(q)}(x)-x\in\OK[x]$, hence integral, and so lies in
some $U_\xi$.
\end{proof}

Orbits of points defined over a finite extension of $K_\p$ approach the orbits of the periodic points of their discs, at a rate controlled by the ramification index.

\begin{lemma}\label{lem:approach}
Let $F\subset\CC_p$ be a finite extension of $K_\p$ with absolute ramification index $e_F$, let
$e\ge e_F$, and put $m_0=\lceil\log_de\rceil\ge0$. Let $x\in F$ be integral with infinite orbit, and
$\eta=\eta_{\bar x}$. Then for all $n\ge0$,
$$v\bigl(f^{(n)}(x)-f^{(n)}(\eta)\bigr)\ \ge\ \min\{d^n/e,\ 1\},$$
and $v(f^{(n)}(x)-f^{(n)}(\eta))\ge n-m_0+1$ for $n\ge m_0$.
\end{lemma}

\begin{proof}
Note $e\ge e_F\ge1$, so $m_0\ge0$. Let $q=q_{\bar x}$. Then $v(f^{(q)}(x)-x)$ is the valuation of a non-zero element
of $F$ (non-zero as $x$ is not periodic), and it is positive; so it is at least $1/e_F\ge1/e$. By
Lemma \ref{lem:periodic}, $g_0=v(x-\eta)\ge1/e$. Put $g_n=v(f^{(n)}(x)-f^{(n)}(\eta))$; by Lemma
\ref{lem:loc}(c), $g_{n+1}\ge\rho(g_n)$. As $\rho$ is increasing, induction gives
$g_n\ge\min\{d^n/e,1\}$: if this holds for $n$, then
$g_{n+1}\ge\min\{d\min\{d^n/e,1\},\ 1\}\ge\min\{d^{n+1}/e,1\}$. Since
$d^{m_0}\ge e$, $g_{m_0}\ge1$. If $g\ge1$ then $dg\ge g+1$,
so $\rho(g)=g+1$; hence $g_n\ge n-m_0+1$ for $n\ge m_0$.
\end{proof}

The following invariants control where the degree can drop.

\begin{definition}\label{def:invariants}
With $\Sfin$ as in the introduction (empty if no finite critical point of $f$ is periodic), put
$$\kappa^*_\p=\sup\{\vartheta(\eta):\ \eta\text{ periodic},\ \bar\eta\notin\Sfin\},\qquad
\nu_\p=\lfloor\kappa^*_\p\rfloor,$$
with $\vartheta$ as in Lemma \ref{lem:loc}(d). For $\alpha\in\Kbar$, $e_\p(\alpha)$ denotes the largest
absolute ramification index of a prime of $K(\alpha)$ above $\p$.
\end{definition}

These invariants are finite.

\begin{lemma}\label{lem:kappa}
One has $1\le\kappa^*_\p<\infty$, so $1\le\nu_\p<\infty$. The set $\Sfin$ is finite and $\bar f(\Sfin)=\Sfin$.
\end{lemma}

\begin{proof}
The superattracting cycles are finitely many, as $f'$ has finitely many zeros, and $f$ maps each of them onto itself;
so $\Sfin$ is finite and $\bar f(\Sfin)=\Sfin$. By Lemma \ref{lem:loc}(d), $\vartheta\ge1$, and
$\vartheta(\eta_\xi)=\vartheta_0$ outside finitely many $\xi$; in particular the set in the definition of $\kappa^*_\p$ is
non-empty. For each of the finitely many remaining $\xi\notin\Sfin$, $\eta_\xi$ is periodic and not in a superattracting
cycle, hence not critical, so $\vartheta(\eta_\xi)<\infty$.
Since the periodic points are exactly the $\eta_\xi$ (Lemma \ref{lem:periodic}), $\kappa^*_\p$ is a maximum of finitely
many finite numbers.
\end{proof}

For $x^d+c$ the invariant $\kappa^*_\p$ can be computed explicitly.

\begin{lemma}\label{lem:kappaxd}
Let $f=x^d+c$ with $c\in\OK$, $d=p^r$, and suppose $0$ is not periodic. Let $q_0\ge1$ be least with
$f^{(q_0)}(0)\in\p$. Then $\kappa^*_\p=r+(d-1)\,v\bigl(f^{(q_0)}(0)\bigr)$.
\end{lemma}

\begin{proof}
Here $f'(\eta)=d\eta^{d-1}$, so $v(f'(\eta_\xi))=r+(d-1)v(\eta_\xi)$, which equals $r$ unless $\xi=0$.
The exact period of $\bar0$ is $q_0$, since $f^{(j)}(0)\in\OK$ lies in the maximal ideal of $\OO_{\CC_p}$
if and only if it lies in $\p$. As $0$ is not periodic, Lemma \ref{lem:periodic} with $x=0$ gives
$v(\eta_0)=v(f^{(q_0)}(0))$.
\end{proof}

We give some values.

\begin{example}\label{ex:kappa}
For $f=x^2+c$ with $c\in\ZZ\setminus\{0,-1\}$ and $p=2$: if $c$ is even then $q_0=1$ and
$\kappa^*_2=1+v_2(c)$; if $c$ is odd then $q_0=2$, $f^{(2)}(0)=c(c+1)$ and $\kappa^*_2=1+v_2(c+1)$. Thus
$\kappa^*_2=2$ for $x^2+1$ and $x^2-2$, and $\kappa^*_2=3$ for $x^2+3$. For $x^p+c$ over $\QQ$ with
$p\nmid c$ and $0$ not periodic, $q_0=p$ and $\kappa^*_p=1+(p-1)v_p(f^{(p)}(0))$; with $p\mid c$, $q_0=1$
and $\kappa^*_p=1+(p-1)v_p(c)$.
\end{example}

\section{Where the degree can drop, and the proof of Theorem \texorpdfstring{A(i)}{A(i)}}\label{sec:drops}

In this section we prove Theorem \ref{thm:droploc}, and then deduce Theorem \ref{thm:main}(i) from it and
Proposition \ref{thm:window}. Throughout, the notation is that of \S\ref{sec:local}; in particular $\Kbar\subset\CC_p$ is fixed, $v(p)=1$, and
$\alpha_n=f^{(n)}(\alpha)$.

\subsection{Outline of the proof}
Suppose that the degree drops at step $n$, that is, $\alpha_n$ has more conjugates than $\alpha_{n+1}$. Then two different
conjugates of $\alpha_n$ have the same image under $f$; call them $x_n\ne y_n$, with $f(x_n)=f(y_n)$. We bound $n$ as
follows.
\begin{enumerate}[label=(\roman*),leftmargin=*]
\item The points $x_n$ and $y_n$ are integral (we shall see that this is automatic), and the reduction $\bar f$ is injective
(Lemma \ref{lem:loc}(a)); so $x_n$ and $y_n$ lie in the same residue disc. Let
$\eta$ be the periodic point of this disc (Lemma \ref{lem:periodic}).
\item Both $x_n$ and $y_n$ are $n$-th iterates of conjugates of $\alpha$, and orbits approach the periodic points: by
Lemma \ref{lem:approach}, for $n\ge m_0$ the points $x_n$ and $y_n$ are within $p^{-(n-m_0+1)}$ of $\eta$, where $m_0=\lceil\log_de_\p(\alpha)\rceil$.
\item Since $f(x_n)=f(y_n)$, the divided difference $(f(x_n)-f(y_n))/(x_n-y_n)$ vanishes. But the divided difference of $f$ at two
points close to $\eta$ is close to $f'(\eta)$. So $f'(\eta)$ is divisible by a high power of $p$, namely
$\vartheta(\eta)=v(f'(\eta))\ge n-m_0+1$.
\end{enumerate}
If $\bar\eta\notin\Sfin$, then $\vartheta(\eta)\le\kappa^*_\p$ by the definition of $\kappa^*_\p$, and we get $n\le m_0+\kappa^*_\p-1$. The
hypothesis $\bar\alpha'\notin\Sfin$ guarantees $\bar\eta\notin\Sfin$, because $\bar\eta=\bar f^n(\bar\alpha')$ and the complement of $\Sfin$ is
stable under $\bar f$.

\subsection{Auxiliary lemmas}

A drop of the degree produces two conjugates of $\alpha$ whose $n$-th iterates are distinct but have the same image under $f$.

\begin{lemma}\label{lem:twoconj}
Let $K\subseteq F\subseteq\Kbar$ be a field, $\alpha\in\Kbar$, and $\alpha'$ a $K$-conjugate of $\alpha$. If $[F(\alpha_{n+1}):F]<[F(\alpha_n):F]$,
then there are $K$-conjugates $x,y$ of $\alpha$, with $x=\alpha'$, such that
$$f^{(n)}(x)\ \ne\ f^{(n)}(y),\qquad f^{(n+1)}(x)\ =\ f^{(n+1)}(y).$$
Moreover, the fields $K_\p(x)$ and $K_\p(y)$ in $\CC_p$ have absolute ramification index at most $e_\p(\alpha)$.
\end{lemma}

\begin{proof}
Since $F(\alpha_{n+1})\subseteq F(\alpha_n)$ and the degrees differ, the inclusion is strict, so by Galois theory there is
$\tau\in\mathrm{Gal}(\Kbar/F(\alpha_{n+1}))$ with $\tau(\alpha_n)\ne\alpha_n$. Fix $\varrho\in\mathrm{Gal}(\Kbar/K)$ with $\varrho(\alpha)=\alpha'$, and put
$x=\varrho(\alpha)=\alpha'$ and $y=\varrho\tau(\alpha)$. As $f$ has coefficients in $K$, which is fixed by $\varrho$ and $\tau$,
$$f^{(n)}(x)=\varrho(\alpha_n)\ne\varrho\tau(\alpha_n)=f^{(n)}(y),\qquad f^{(n+1)}(y)=\varrho\tau(\alpha_{n+1})=\varrho(\alpha_{n+1})=f^{(n+1)}(x),$$
using $\tau(\alpha_{n+1})=\alpha_{n+1}$. Finally, the field $K_\p(x)\subset\CC_p$ is the completion of $K(\alpha)$ at the prime above
$\p$ induced by the embedding $\varrho|_{K(\alpha)}\colon K(\alpha)\to\CC_p$, so its absolute ramification index is that of a prime of $K(\alpha)$
above $\p$, hence at most $e_\p(\alpha)$; likewise for $y$, with the embedding $\varrho\tau|_{K(\alpha)}$.
\end{proof}

These two iterates lie in the same residue disc and are close to its periodic point.

\begin{lemma}\label{lem:samedisc}
Let $\alpha\in\Kbar$ have infinite $f$-orbit, and let $x,y\in\CC_p$ be $K$-conjugates of $\alpha$ with $x$ integral, $\bar x\notin\Sfin$, and suppose $x_n=f^{(n)}(x)$ and
$y_n=f^{(n)}(y)$ satisfy $f(x_n)=f(y_n)$. Then $y$ is integral, $x_n$ and $y_n$ lie in the same residue disc $U_\xi$ with
$\xi\notin\Sfin$, and $f^{(n)}(\eta_{\bar x})=\eta_\xi=f^{(n)}(\eta_{\bar y})$. If moreover the fields $K_\p(x)$ and $K_\p(y)$ have
absolute ramification index at most $e$, and $m_0=\lceil\log_de\rceil$, then for $n\ge m_0$
$$v(x_n-\eta_\xi)\ \ge\ n-m_0+1,\qquad v(y_n-\eta_\xi)\ \ge\ n-m_0+1.$$
\end{lemma}

\begin{proof}
Since $f$ is monic with integral coefficients, $x_n$ and $f(x_n)$ are integral. If $y_n$ were not
integral, then by Lemma \ref{lem:loc}(e) $f(y_n)$ would not be integral either, contradicting $f(y_n)=f(x_n)$. So $y_n$ is
integral, and then so is $y$: if $y$ were not integral, none of its iterates would be, again by Lemma \ref{lem:loc}(e).

Reducing $f(x_n)=f(y_n)$ gives $\bar f(\bar x_n)=\bar f(\bar y_n)$, and $\bar f$ is injective by Lemma
\ref{lem:loc}(a); so $\bar x_n=\bar y_n=:\xi$. Since $\xi=\bar f^n(\bar x)$, $\bar x\notin\Sfin$ and the complement of $\Sfin$ is stable under the
bijection $\bar f$ (Lemma \ref{lem:kappa}), we have $\xi\notin\Sfin$.

By Lemma \ref{lem:periodic}, $f(\eta_\zeta)=\eta_{\bar f(\zeta)}$ for every $\zeta\in\overline{\FF}_p$; iterating,
$f^{(n)}(\eta_{\bar x})=\eta_{\bar f^n(\bar x)}=\eta_\xi$, and likewise for $y$, since $\bar f^n(\bar y)=\bar y_n=\xi$.

Apply Lemma \ref{lem:approach} to $x$, with $F=K_\p(x)$, which is a finite extension of
$K_\p$ with ramification index at most $e$; as $x$ has infinite orbit (it is conjugate to $\alpha$), it gives
$v(x_n-f^{(n)}(\eta_{\bar x}))\ge n-m_0+1$ for $n\ge m_0$, and $f^{(n)}(\eta_{\bar x})=\eta_\xi$. The same argument applies to $y$.
\end{proof}

Finally, a divided difference of $f$ near $\eta$ is close to $f'(\eta)$.

\begin{lemma}\label{lem:divdiff}
Let $\eta\in\OO_{\CC_p}$, and let $s\ne t$ in $\OO_{\CC_p}$ satisfy $f(\eta+s)=f(\eta+t)$. Then
$$\vartheta(\eta)=v\bigl(f'(\eta)\bigr)\ \ge\ \min\{v(s),v(t)\}.$$
\end{lemma}

\begin{proof}
Let $A(u)=f(\eta+u)-f(\eta)$, a polynomial in the variable $u$ with coefficients in $\OO_{\CC_p}$, as $\eta$ is integral, with $A(0)=0$ and
$A'(0)=f'(\eta)$. The divided difference
$$G(u_1,u_2)=\frac{A(u_1)-A(u_2)}{u_1-u_2}$$
is a polynomial with coefficients in $\OO_{\CC_p}$: if $A(u)=\sum_iA_iu^i$, then
$G=\sum_iA_i(u_1^{i-1}+u_1^{i-2}u_2+\dots+u_2^{i-1})$. Its constant term is $A_1=A'(0)$, so $G(0,0)=f'(\eta)$, and
$G(u_1,u_2)-G(0,0)$ is a polynomial without constant term and with integral coefficients. Hence
$$v\bigl(G(s,t)-G(0,0)\bigr)\ \ge\ \min\{v(s),v(t)\}.$$
On the other hand $G(s,t)=(A(s)-A(t))/(s-t)=0$, because $A(s)=f(\eta+s)-f(\eta)=f(\eta+t)-f(\eta)=A(t)$ and $s\ne t$. So
$v(f'(\eta))=v(G(0,0))=v(G(0,0)-G(s,t))\ge\min\{v(s),v(t)\}$.
\end{proof}

\subsection{Proof of Theorem \ref{thm:droploc}}

\begin{proof}[Proof of Theorem \ref{thm:droploc}]
Let $e=e_\p(\alpha)$ and $m_0=\lceil\log_de\rceil$, and suppose that $[F(\alpha_{n+1}):F]<[F(\alpha_n):F]$. We must show
$n\le m_0+\nu_\p-1$. If $n<m_0$ this holds, because $\nu_\p\ge1$ (Lemma \ref{lem:kappa}); so assume $n\ge m_0$.

By Lemma \ref{lem:twoconj}, applied with the conjugate $\alpha'$ of the hypothesis, there are $K$-conjugates $x=\alpha'$ and $y$ of
$\alpha$ with $x_n\ne y_n$ and $f(x_n)=f(y_n)$, where $x_n=f^{(n)}(x)$ and $y_n=f^{(n)}(y)$, and the fields $K_\p(x)$, $K_\p(y)$ have
ramification index at most $e$. As $\alpha$ is not preperiodic, it has infinite orbit. By hypothesis $x$ is integral and $\bar x\notin\Sfin$. So Lemma \ref{lem:samedisc} applies: there is
$\xi\notin\Sfin$ with $x_n,y_n\in U_\xi$ and, writing $\eta=\eta_\xi$,
$$v(x_n-\eta)\ \ge\ n-m_0+1,\qquad v(y_n-\eta)\ \ge\ n-m_0+1.$$
Put $s=x_n-\eta$ and $t=y_n-\eta$. Then $s\ne t$, and $f(\eta+s)=f(x_n)=f(y_n)=f(\eta+t)$, so Lemma \ref{lem:divdiff} gives
$$\vartheta(\eta)\ \ge\ \min\{v(s),v(t)\}\ \ge\ n-m_0+1.$$
Now $\eta$ is periodic and $\bar\eta=\xi\notin\Sfin$, so $\vartheta(\eta)\le\kappa^*_\p$ by Definition \ref{def:invariants}. Hence $n-m_0+1\le\kappa^*_\p$,
and since $n-m_0+1$ is an integer, $n-m_0+1\le\lfloor\kappa^*_\p\rfloor=\nu_\p$, that is, $n\le m_0+\nu_\p-1$.
\end{proof}

\subsection{Proof of Theorem \ref{thm:main}(i)}

\begin{proof}[Proof of Theorem \ref{thm:main}\textup{(i)}]
Let $\alpha\in\Kbar$ be not preperiodic, with $D=[K(\alpha):K]\ge2$.

If $\alpha$ is not integral over $\OK$, Lemma \ref{lem:nonintegral} gives
$\hh_f(\alpha)\ge\log2/([K:\QQ]D)$. This is larger than the right side of (i) as soon as
$c_1\le\log2/[K:\QQ]$, since $e_\p(\alpha)\ge1$ and the exponential factor is at most $1$.

Let $\alpha$ be integral over $\OK$. Then every $K$-conjugate $\alpha'$ is integral, so
$v(\alpha')\ge0$; and $\Sfin$ is empty, because no finite critical point is periodic. So the hypotheses of Theorem
\ref{thm:droploc} hold, and applying it with $F=K$ we find that $D_{n+1}=D_n$ for every
$$n\ \ge\ m:=\lceil\log_de_\p(\alpha)\rceil+\nu_\p.$$

Let $\ell=\log3D$ and $L=\lceil a\sqrt{q\ell}\,\rceil$, as in Proposition \ref{thm:window}.
The degrees $D_n$ are constant for $n\ge m$, in particular for $m\le n\le m+n_0+L$, and Proposition
\ref{thm:window} applies with $j=m$:
$$\hh_f(\alpha)\ \ge\ \frac{c_3}Dd^{-a\sqrt{q\ell}}\prod_{n<m}\frac{k_n}d.$$
Each factor satisfies $k_n/d\ge1/d$, so the product is at least $d^{-m}$, and
$$d^{\,m}\ \le\ d^{\,\log_de_\p(\alpha)+1+\nu_\p}\ =\ d^{\,1+\nu_\p}\,e_\p(\alpha).$$
Hence $\hh_f(\alpha)\ge c_3d^{-1-\nu_\p}\,e_\p(\alpha)^{-1}D^{-1}\exp(-a\sqrt q\log d\,\sqrt{\log3D})$, which is (i)
for every $c_1\le\min\{c_3d^{-1-\nu_\p},\log2/[K:\QQ]\}$, with $c_2=a\sqrt q\log d$. The constant $c_1$ of Theorem \ref{thm:main} is the
minimum of this bound and of those required in the proof of (ii) in \S\ref{sec:clusters}.
\end{proof}

\section{Ramified clusters and the proof of Theorem \texorpdfstring{A(ii)}{A(ii)}}\label{sec:clusters}

In this section we prove part (ii) of Theorem \ref{thm:main}, in a form that also applies when some finite critical point is
periodic. Throughout, $f\in\OK[x]$ is monic of Frobenius shape at $p$, the notation is that of \S\ref{sec:local}, and
$\Lambda(N)$ is as in \S\ref{sec:proofA}. For a $K$-conjugate $x\in\Kbar\subset\CC_p$ of $\alpha$ let $e_x$ be the ramification
index of $K_\p(x)$ over $\QQ_p$; thus $e_x\le e_\p(\alpha)$.

\subsection{A cluster of conjugates}

We first determine how the conjugates of $\alpha$ distribute over those of $\alpha_n$.

\begin{lemma}\label{lem:fibres}
Let $\alpha\in\Kbar$ and $n\ge0$. The map $x\mapsto f^{(n)}(x)$ sends the set of $K$-conjugates of $\alpha$ onto the set of
$K$-conjugates of $\alpha_n$, and each fibre has exactly $D/D_n$ elements. In particular, the image of a set $C$ of
$K$-conjugates of $\alpha$ has at least $|C|D_n/D$ elements.
\end{lemma}

\begin{proof}
Let $G=\mathrm{Gal}(\Kbar/K)$. Since $f$ has coefficients in $K$, $f^{(n)}(g\alpha)=g\alpha_n$ for $g\in G$, so the map is
$G$-equivariant; it is onto, because $G$ acts transitively on the conjugates of $\alpha_n$. If $g\alpha_n=y$, then $g$
maps the fibre over $\alpha_n$ bijectively onto the fibre over $y$, so all fibres have the same size, namely $D/D_n$. A
set $C$ meets each fibre in at most $D/D_n$ points, so its image has at least $|C|/(D/D_n)$ elements.
\end{proof}

Each conjugate $x$ of $\alpha$ lies in a cluster of $e_x/e(\p|p)$ conjugates whose orbits remain $\p$-adically close.

\begin{lemma}\label{lem:cluster}
Let $\alpha\in\Kbar$ be integral over $\OK$ with infinite orbit, and let $x$ be a $K$-conjugate of $\alpha$. There is a set $C(x)$ of
$K$-conjugates of $\alpha$ containing $x$, with
$$|C(x)|\ =\ \frac{e_x}{e(\p|p)},$$
such that for all $x'\ne x''$ in $C(x)$ and all $n\ge0$,
$$v\bigl(f^{(n)}(x')-f^{(n)}(x'')\bigr)\ \ge\ \min\{d^n/e_x,\ 1\}.$$
\end{lemma}

\begin{proof}
Let $E=K_\p(x)\subset\CC_p$, and let $T\subseteq E$ be the maximal unramified extension of $K_\p$ contained in $E$. Then $E=T(x)$,
$[E:T]=e_x/e(\p|p)$, and $E$ and $T$ have the same residue field. Let $C(x)$ be the set of roots in $\CC_p$ of the minimal polynomial
of $x$ over $T$. It has $[E:T]$ elements, as we are in characteristic $0$, and since this minimal polynomial divides that of $\alpha$
over $K$, $C(x)$ consists of $K$-conjugates of $\alpha$, all of them integral.

Let $x'\in C(x)$, and let $\varphi\colon T(x)\to T(x')$ be the $T$-isomorphism with $\varphi(x)=x'$. As $T$ is complete, $v$ has a
unique extension from $T$ to each finite extension, so $v\circ\varphi=v$ on $T(x)$. In particular $K_\p(x')=T(x')$ has
ramification index $e_x$. Since $E$ and $T$ have the same residue field, there is $t\in\OO_T$ with $v(x-t)>0$; then
$v(x'-t)=v(\varphi(x-t))=v(x-t)>0$, so $\bar x'=\bar x$. Hence all points of $C(x)$ lie in the residue disc of $\xi=\bar x$; let
$\eta=\eta_\xi$.

Each $x'\in C(x)$ has infinite orbit, being conjugate to $\alpha$. So Lemma \ref{lem:approach}, applied with $F=K_\p(x')$, gives
$v(f^{(n)}(x')-f^{(n)}(\eta))\ge\min\{d^n/e_x,1\}$ for every $x'\in C(x)$, with the same $\eta$, and the ultrametric inequality
gives the claim.
\end{proof}

\subsection{Drops force ramification of many conjugates}

A drop of the degree at step $n$ forces at least half of the conjugates reducing outside $\Sfin$ to be highly ramified above $\p$.

\begin{lemma}\label{lem:half}
Let $\alpha$ be integral over $\OK$ and not preperiodic, suppose that the degree drops at step $n$, and let $X$ be the set of the
$K$-conjugates of $\alpha$ that reduce outside $\Sfin$. Then the set
$$H_n\ =\ \{x\in X\ :\ e_x>d^{\,n-\nu_\p}\}$$
has at least $|X|/2$ elements.
\end{lemma}

\begin{proof}
By Lemma \ref{lem:fibres}, $f$ maps the conjugates of $\alpha_n$ onto those of $\alpha_{n+1}$, with fibres of $k_n\ge2$ elements. The
points of a fibre have the same image, so they lie in the same residue disc, as $\bar f$ is injective (Lemma \ref{lem:loc}(a)); and
since $\bar f(\Sfin)=\Sfin$ and $\bar f$ is a bijection (Lemma \ref{lem:kappa}), a conjugate $x$ of $\alpha$ reduces outside $\Sfin$ if and
only if $f^{(n)}(x)$ does. Call a conjugate $z$ of $\alpha_n$ \emph{bad} if some conjugate $x$ of $\alpha$ with $f^{(n)}(x)=z$
satisfies $e_x\le d^{\,n-\nu_\p}$.

We claim that a fibre $\Phi$ outside $\Sfin$ contains at most one bad element. Suppose $z\ne z'$ in $\Phi$ are bad, with $x,x'$ as in
the definition; as $e_x\ge1$, this forces $n\ge\nu_\p$. All these points are integral. Let $\eta$ be the periodic point of the residue
disc of $z$ and $z'$; then $\bar\eta\notin\Sfin$ and $f^{(n)}(\eta_{\bar x})=\eta=f^{(n)}(\eta_{\bar x'})$ (Lemma \ref{lem:periodic}). Put
$m=\lceil\log_de_x\rceil$; since $\log_de_x\le n-\nu_\p$ and $n-\nu_\p$ is an integer, $m\le n-\nu_\p<n$, and Lemma
\ref{lem:approach}, applied with $F=K_\p(x)$, gives $v(z-\eta)\ge n-m+1\ge\nu_\p+1$; likewise $v(z'-\eta)\ge\nu_\p+1$. Lemma
\ref{lem:divdiff}, with $s=z-\eta\ne t=z'-\eta$ and $f(\eta+s)=f(\eta+t)$, gives $\vartheta(\eta)\ge\nu_\p+1$, whereas
$\vartheta(\eta)\le\kappa^*_\p<\nu_\p+1$ by Definition \ref{def:invariants}. This contradiction proves the claim.

The conjugates of $\alpha_n$ outside $\Sfin$ form $|X|D_n/(Dk_n)$ fibres, so at least $(1-1/k_n)|X|D_n/D\ge|X|D_n/(2D)$ of them are not
bad, and all $D/D_n$ conjugates of $\alpha$ above each of them lie in $H_n$. Hence $|H_n|\ge|X|/2$.
\end{proof}

If $X$ contains all conjugates, Theorem \ref{thm:droploc} for $F=K$ is the statement that $H_n$ is not empty.

\subsection{The discriminant at the last drop}

At a drop $m$, these clusters make the discriminant of the minimal polynomial of $\alpha_m$ divisible by a large power of $p$, which gives the following bound.

\begin{proposition}\label{prop:disc}
Put $c_8=\log p/(4[K:\QQ]\,e(\p|p)\,d^{\,\nu_\p})$, and let $c_9\ge1$ be such that $\Lambda(N)+\log p\le c_9\log3N$ for all
$N\ge1$. Let $\alpha$ be integral over $\OK$ and not preperiodic, let $\rho$ be the proportion of the $K$-conjugates of $\alpha$
that reduce outside $\Sfin$, and suppose that the degree drops at step $m$. Then
$$\hh_f(\alpha)\ \ge\ c_8\,\frac{\rho}{D}\ -\ c_9\,\frac{\log3D_m}{D_m\,d^{\,m}}.$$
\end{proposition}

\begin{proof}
Let $P\in\OK[x]$ be the minimal polynomial of $\alpha_m$ over $K$, monic of degree $D_m$, and $\Delta=\disc(P)\in\OK$, which is
non-zero.

We first bound $\Delta$ at $\p$. Let $H=H_m$ be as in Lemma \ref{lem:half}, so that $|H|\ge\rho D/2$, and let $x\in H$. By Lemma
\ref{lem:cluster}, the set $f^{(m)}(C(x))$ of roots of $P$ contains $y=f^{(m)}(x)$, its points are pairwise at valuation
$\ge\delta_x=\min\{d^m/e_x,1\}$, and by Lemma \ref{lem:fibres} it has at least $|C(x)|D_m/D$ elements. So $y$ has at least
$|C(x)|D_m/D-1$ other roots $y'$ of $P$ with $v(y-y')\ge\delta_x$, and
$$\sum_{y'\ne y}v(y-y')\ \ge\ \Bigl(\frac{|C(x)|D_m}D-1\Bigr)\delta_x\ \ge\ \frac{\min\{d^m,e_x\}}{e(\p|p)}\cdot\frac{D_m}D-1\ \ge\
\frac{d^{\,m-\nu_\p}D_m}{e(\p|p)\,D}-1,$$
using $\delta_x\le1$ and $e_x>d^{\,m-\nu_\p}$. Since $\Delta=\pm\prod_{y\ne y'}(y-y')$ over ordered pairs of roots, all integral, $v(\Delta)$ is
at least the sum of these quantities over the roots $y$ of the form $f^{(m)}(x)$ with $x\in H$. Each root has $D/D_m$ conjugates of
$\alpha$ above it, so there are at least $|H|D_m/D\ge\rho D_m/2$ such roots, and, as $v(\Delta)\ge0$,
$$v(\Delta)\ \ge\ \frac{\rho D_m}2\Bigl(\frac{d^{\,m-\nu_\p}D_m}{e(\p|p)\,D}-1\Bigr).$$
The norm $\Nm_{K/\QQ}(\Delta)$ is a non-zero rational integer, and $v_p(\Nm_{K/\QQ}\Delta)\ge e(\p|p)f(\p|p)\,v(\Delta)\ge v(\Delta)$,
because $v_\p=e(\p|p)\,v$ on $K$ and all $v_{\mathfrak q}(\Delta)\ge0$. Therefore $\log|\Nm_{K/\QQ}\Delta|\ge v(\Delta)\log p$.

Next we bound $\Delta$ at the archimedean places. Fix $\sigma\colon K\to\CC$, let $Z_\sigma$ be the set of roots of $\sigma P$, and let $\Delta_\sigma$ be the
determinant of the matrix $(\sigma B_k(z))_{z\in Z_\sigma,\,0\le k<D_m}$. As in Lemma \ref{lem:detint}, $\Delta_\sigma^2=\pm\sigma(\Delta)$; and
exactly as in the proof of Lemma \ref{lem:detup}, with the single level $\alpha_m$ and $N=D_m$, using \eqref{eq:hsum},
$$\tfrac12\log\bigl|\Nm_{K/\QQ}\Delta\bigr|\ =\ \sum_\sigma\log|\Delta_\sigma|\ \le\ [K:\QQ]\,D_m\Lambda(D_m)+[K:\QQ]\,D_m^2d^{\,m}\,\hh_f(\alpha).$$

Combining the two bounds and dividing by $[K:\QQ]D_m^2d^{\,m}$,
$$\hh_f(\alpha)\ \ge\ \frac{\rho\log p}{4[K:\QQ]\,e(\p|p)\,d^{\,\nu_\p}D}-\frac{\Lambda(D_m)+\log p}{D_m\,d^{\,m}},$$
which is the claim.
\end{proof}

Combining Proposition \ref{prop:disc} with the bound for the iterates after the last drop, whose degrees are constant, we obtain the following bound for points with a proportion $\rho$ of conjugates reducing outside $\Sfin$.

\begin{proposition}\label{prop:outside}
There are $c_{10},c_2>0$, depending only on $K$, $f$ and $\p$, such that every $\alpha$ integral over $\OK$ and not preperiodic, with
$D\ge2$ and with a proportion $\rho>0$ of its $K$-conjugates reducing outside $\Sfin$, satisfies
$$\hh_f(\alpha)\ \ge\ \frac{c_{10}\,\rho}{D\log3D}\exp\bigl(-c_2\sqrt{\log3D}\bigr).$$
\end{proposition}

\begin{proof}
Let $\ell=\log3D$, let $a,q,c_3$ be as in Proposition \ref{prop:window}, $c_2=a\sqrt q\log d$, and $H_D=d^{-a\sqrt{q\ell}}=\exp(-c_2\sqrt\ell)$.
If the degree never drops, Proposition \ref{thm:window} with $j=0$ gives $\hh_f(\alpha)\ge c_3H_D/D$. Otherwise, as $(D_n)$ is a
non-increasing sequence of positive integers, there is a last step $m$ at which the degree drops, and Proposition \ref{thm:window}
with $j=m+1$ gives
$$\hh_f(\alpha)\ \ge\ \frac{c_3H_D}D\prod_{n\le m}\frac{k_n}d\ =\ \frac{c_3H_D}{D_{m+1}\,d^{\,m+1}}\ \ge\ \frac{c_3}d\cdot\frac{H_D}{D_m\,d^{\,m}}.$$
If $c_8\rho/D\ge2c_9\log(3D_m)/(D_md^{\,m})$, Proposition \ref{prop:disc} gives $\hh_f(\alpha)\ge c_8\rho/(2D)$. Otherwise
$D_md^{\,m}<(2c_9/c_8)\rho^{-1}D\log3D$, and the bound for iterates of constant degree gives $\hh_f(\alpha)\ge\frac{c_3c_8}{2c_9d}\cdot\frac{\rho H_D}{D\log3D}$. The claim
holds with $c_{10}=\min\{c_3,\,c_8/2,\,c_3c_8/(2c_9d)\}$, as $\rho\le1$ and $H_D\le1\le\log3D$.
\end{proof}

\begin{proof}[Proof of Theorem \ref{thm:main}\textup{(ii)}]
If $\alpha$ is not integral over $\OK$, Lemma \ref{lem:nonintegral} gives $\hh_f(\alpha)\ge\log2/([K:\QQ]D)$. Otherwise, since no finite
critical point is periodic, $\Sfin$ is empty and $\rho=1$, and Proposition \ref{prop:outside} applies. Since $\sqrt y\ge\log y$ for $y>0$,
we have $\log3D\le\exp(\sqrt{\log3D})$, so Proposition \ref{prop:outside} gives $\hh_f(\alpha)\ge c_{10}D^{-1}\exp\bigl(-(c_2+1)\sqrt{\log3D}\bigr)$.
So (ii) holds with $c_1\le\min\{c_{10},\log2/[K:\QQ]\}$.
\end{proof}

\section{Superattracting classes}\label{sec:attract}

In this section we treat the integral points all of whose conjugates reduce into $\Sfin$, using that $f$ contracts $\p$-adically near superattracting cycles, and we prove Theorems \ref{thm:attract} and \ref{thm:all}. Here $f\in\OK[x]$ is a monic polynomial of Frobenius shape at $p$ with a superattracting cycle in $\Kbar$, and
$m_Z$, $Q$, $M$, $m_f$, $\omega_f$ and $C_f$ are as in Definition \ref{def:omega}. By Lemma \ref{lem:periodic},
distinct periodic points have distinct reductions, so $|\Sfin|=M$, and $\bar\alpha'\in\Sfin$ means that $\alpha'$ lies in the
residue disc of a point of a superattracting cycle.

Near a superattracting cycle, $f$ contracts $\p$-adic distances to the cycle.

\begin{lemma}\label{lem:contraction}
Let $Z$ be a superattracting cycle of exact period $q$, $\gamma\in Z$, and $x\in\OO_{\CC_p}$ with $\bar x=\bar\gamma$. Put
$s_j=v\bigl(f^{(j)}(x)-f^{(j)}(\gamma)\bigr)\in(0,\infty]$. Then for all $j\ge0$
$$s_{j+1}\ \ge\ \min\{d\,s_j,\ 1\},\qquad s_{j+1}\ \ge\ s_j,\qquad\text{and}\qquad s_{j+q}\ \ge\ m_Zs_j.$$
\end{lemma}

\begin{proof}
For $y,z\in\OO_{\CC_p}$, $f(y)-f(z)=(y-z)\,u(y,z)$ with $u\in\OO_{\CC_p}$, as $f$ has integral coefficients; hence
$v(f(y)-f(z))\ge v(y-z)$, which gives $s_{j+1}\ge s_j$. By Lemma \ref{lem:loc}(c), $s_{j+1}\ge\min\{ds_j,s_j+1\}\ge\min\{ds_j,1\}$. Now
let $c\in Z$ be critical with $k=\deg_cf$. As $c$ is integral, $f(c+t)-f(c)=\sum_{i=k}^db_it^i$ with $b_i\in\OO_{\CC_p}$, so
$v(f(c+t)-f(c))\ge k\,v(t)$ for $v(t)\ge0$. Among $q$ consecutive steps the orbit of $\gamma$ passes through every critical point $c$ of
$Z$; at such a step $s$ is multiplied by at least $\deg_cf$, and at the other steps it does not decrease. Hence $s_{j+q}\ge m_Zs_j$.
\end{proof}

A single resultant between $\alpha_n$ and $\alpha_{n+Q}$ then gives a bound in terms of the ramification, for points with a proportion $\sigma$ of conjugates reducing into $\Sfin$.

\begin{proposition}\label{prop:attractmix}
Let $\alpha\in\Kbar$ be integral over $\OK$ and not preperiodic, let $e=e_\p(\alpha)$, and suppose that a proportion $\sigma>0$ of its
$K$-conjugates reduce into $\Sfin$. With $A$ and $B$ as in Theorem \ref{thm:attract},
$$\hh_f(\alpha)\ \ge\ \frac B{2d^{2Q+1}}\min\Bigl\{1,\Bigl(\frac{\sigma^2}{2AB}\Bigr)^{\omega_f}\Bigr\}\cdot\frac1e.$$
\end{proposition}

\begin{proof}
Let $\alpha^{(1)},\dots,\alpha^{(D)}$ be the $K$-conjugates of $\alpha$, numbered so that $\alpha^{(i)}$ reduces into $\Sfin$ exactly for
$i\le\sigma D$, and for $j\ge0$ let
$$\chi_j(X)=\prod_{i=1}^D\bigl(X-f^{(j)}(\alpha^{(i)})\bigr)\in\OK[X],$$
the characteristic polynomial of $f^{(j)}(\alpha)$ over $K$, raised to the power $D/[K(f^{(j)}(\alpha)):K]$.

By Lemma \ref{lem:periodic}, for each $i\le\sigma D$ there is a unique point $\gamma_i$ of a superattracting cycle $Z_i$ with
$\bar\alpha^{(i)}=\bar\gamma_i$. The field $F_i=K_\p(\alpha^{(i)})$ has ramification index $e_i\le e$ over $\QQ_p$, and $\gamma_i\in F_i$ by
Lemma \ref{lem:periodic}. As $\alpha$ is not preperiodic, $\alpha^{(i)}-\gamma_i$ is a non-zero element of $F_i$ of positive valuation, so
$s_0^{(i)}=v(\alpha^{(i)}-\gamma_i)\ge1/e$, where $s_j^{(i)}=v\bigl(f^{(j)}(\alpha^{(i)})-f^{(j)}(\gamma_i)\bigr)$. Let $j_0=\lceil\log_de\rceil$, so
that $e\le d^{\,j_0}<de$. By Lemma \ref{lem:contraction}, $s_j^{(i)}\ge\min\{d^j/e,1\}$ by induction, so $s_{j_0}^{(i)}\ge1$; and since the period
of $Z_i$ divides $Q$ and $m_{Z_i}\ge m_f$,
\begin{equation}\label{eq:growth}s_{j_0+tQ}^{(i)}\ \ge\ m_f^{\,t}\qquad(t\ge0).\end{equation}

Now fix $n\ge0$, put $m=n+Q$, and consider the resultant
$R=\mathrm{Res}(\chi_n,\chi_m)=\prod_{i,i'}\bigl(f^{(n)}(\alpha^{(i)})-f^{(m)}(\alpha^{(i')})\bigr)$. It lies in $\OK$, and $R\ne0$:
otherwise $f^{(n)}(\alpha)$ and $f^{(m)}(\alpha)$ would be conjugate over $K$, so $d^n\hh_f(\alpha)=d^m\hh_f(\alpha)$,
$\hh_f(\alpha)=0$, and $\alpha$ would be preperiodic \cite[Thm.~3.22]{Sil07}.

We bound $R$ from below at $\p$. For $i\le\sigma D$ put $\delta_i=f^{(n)}(\gamma_i)$, one of the $M$ points of the superattracting
cycles. If $\delta_i=\delta_{i'}$, then $f^{(m)}(\gamma_{i'})=f^{(n)}(\gamma_{i'})=\delta_{i'}$, because the period of $Z_{i'}$
divides $Q$; so by the ultrametric inequality,
$$v\bigl(f^{(n)}(\alpha^{(i)})-f^{(m)}(\alpha^{(i')})\bigr)\ \ge\ \min\bigl(s_n^{(i)},s_m^{(i')}\bigr)\ \ge\ s_n,\qquad s_n=\min_{i\le\sigma D}s_n^{(i)}.$$
All other factors of $R$ are integral. If $\mu_\delta$ is the number of $i\le\sigma D$ with $\delta_i=\delta$, the number of ordered pairs
$(i,i')$ with $i,i'\le\sigma D$ and $\delta_i=\delta_{i'}$ is $\sum_\delta \mu_\delta^2\ge\sigma^2D^2/M$ by the Cauchy--Schwarz inequality.
Hence $v(R)\ge(\sigma^2D^2/M)\,s_n$, and since $v_\p=e(\p|p)\,v$ on $K$ and $R\in\OK$,
$$\log\bigl|\Nm_{K/\QQ}(R)\bigr|\ \ge\ v_\p(R)\log\Nm\p=[K_\p:\QQ_p]\log p\cdot v(R)\ \ge\
[K_\p:\QQ_p]\log p\cdot\frac{\sigma^2D^2}M\,s_n.$$

We bound $R$ from above at the archimedean places. Let $\tau$ run over the embeddings of $K$ into $\CC$, each extended to
$\Kbar$. For complex $y,z$ we have $|y-z|\le2\max(1,|y|)\max(1,|z|)$. For an algebraic integer $\beta\in K(\alpha)$, the sum
of $\log\max(1,|\tau\beta^{(i)}|)$ over all $\tau$ and all $K$-conjugates $\beta^{(i)}$ indexed as above is
$[K(\alpha):\QQ]\,h(\beta)=[K:\QQ]\,D\,h(\beta)$. Hence
\begin{align*}
\log\bigl|\Nm_{K/\QQ}(R)\bigr|=\sum_\tau\log|\tau R|
&\le\sum_\tau\sum_{i,i'}\Bigl(\log2+\log^+\bigl|\tau f^{(n)}(\alpha^{(i)})\bigr|+\log^+\bigl|\tau f^{(m)}(\alpha^{(i')})\bigr|\Bigr)\\
&=[K:\QQ]\,D^2\Bigl(\log2+h\bigl(f^{(n)}(\alpha)\bigr)+h\bigl(f^{(m)}(\alpha)\bigr)\Bigr).
\end{align*}
By the definition of $C_f$ and $\hh_f(f^{(j)}(\alpha))=d^j\hh_f(\alpha)$, the bracket is at most $B+2d^{n+Q}\hh_f(\alpha)$.

Comparing the two bounds and cancelling the factor $D^2$, we obtain
\begin{equation}\label{eq:keyattract}\sigma^2s_n\ \le\ A\,\bigl(B+2d^{n+Q}\hh_f(\alpha)\bigr)\qquad\text{for every }n\ge0.\end{equation}
Let $t\ge0$ be the least integer with $m_f^t>2AB/\sigma^2$, and take $n=j_0+tQ$. By \eqref{eq:growth}, $\sigma^2s_n>2AB$, so \eqref{eq:keyattract}
gives $2Ad^{n+Q}\hh_f(\alpha)>AB$, that is, $\hh_f(\alpha)>B/(2d^{\,j_0+(t+1)Q})$. If $t\ge1$, minimality gives $m_f^{t-1}\le2AB/\sigma^2$, so
$d^{Q(t-1)}\le(2AB/\sigma^2)^{\omega_f}$; in all cases $d^{Qt}\le d^Q\max\{1,(2AB/\sigma^2)^{\omega_f}\}$. With $d^{\,j_0}<de$ this gives the claim.
\end{proof}

\begin{proof}[Proof of Theorem \ref{thm:attract}]
This is Proposition \ref{prop:attractmix} with $\sigma=1$.
\end{proof}

\begin{proof}[Proof of Theorem \ref{thm:all}]
If no finite critical point is periodic, this is Theorem \ref{thm:main}(ii). Otherwise let $\alpha$ be non-preperiodic with $D\ge2$.
If $\alpha$ is not integral over $\OK$, then $\hh_f(\alpha)\ge\log2/([K:\QQ]D)$ by Lemma \ref{lem:nonintegral}. Let $\alpha$ be integral, let
$\rho$ be the proportion of its $K$-conjugates reducing outside $\Sfin$, and $\sigma=1-\rho$. If $\rho\ge1/2$, Proposition
\ref{prop:outside} gives $\hh_f(\alpha)\ge c_{10}\exp(-c_2\sqrt{\log3D})/(2D\log3D)\ge c_{10}\exp(-(c_2+1)\sqrt{\log3D})/(2D)$, since
$\log3D\le\exp(\sqrt{\log3D})$. If $\sigma\ge1/2$, Proposition \ref{prop:attractmix}
gives $\hh_f(\alpha)\ge c/e_\p(\alpha)\ge c/(e(\p|p)D)$, with $c=\frac B{2d^{2Q+1}}\min\{1,(8AB)^{-\omega_f}\}$. In all cases the claim follows, with $c_2+1$ in place of $c_2$.
\end{proof}

\begin{remark}\label{rem:contraction}
Using only $s_{j+q}\ge m_Zs_j$ would give the exponent $\omega_f$ of $e_\p(\alpha)$. By the first inequality of Lemma
\ref{lem:contraction}, $f$ multiplies the valuation of the distance to the cycle by $d$ until it reaches $1$, which takes
$j_0=\lceil\log_de_\p(\alpha)\rceil$ steps, and after that a bounded number of periods suffices. For $f=x^d$ this is the elementary bound
$h(\alpha)\gg1/e_\p(\alpha)$ for integers divisible by every prime above $\p$. The exponent of $e_\p(\alpha)$ is optimal: $\alpha=p^{1/e}$ has
$h(\alpha)=\log p/e$, and for $f=x^2-1$ the points of $f^{-N}(2)$ are totally ramified above $2$ and have $\hh_f(\alpha)=\hh_f(2)/e_2(\alpha)$.
\end{remark}

\section{Examples and explicit constants}\label{sec:examples}

In this section we list polynomials covered by Theorem \ref{thm:main}, and make the constants explicit for $x^2\pm1$.

\subsection{Polynomials covered by the uniform bound}

Theorem \ref{thm:main} applies in particular to the following polynomials.

\begin{corollary}\label{cor:examples}
Theorem \ref{thm:main} applies, with the stated value of $\kappa^*_\p$, to:
\begin{enumerate}[label=\rm(\alph*)]
\item $x^2+c$ with $c\in\ZZ\setminus\{0,-1\}$, $\p=2$: $\kappa^*_2=1+v_2(c)$ for $c$ even,
$\kappa^*_2=1+v_2(c+1)$ for $c$ odd. This includes $x^2+1$, $x^2-2$ and every $c$ for which $x^2+c$ is
not postcritically finite.
\item $x^d+c$ with $d=p^r$, $c\in\OK$ and $0$ not periodic, with $\kappa^*_\p$ as in Lemma
\ref{lem:kappaxd}; e.g.\ $x^2+i$ over $\QQ(i)$, where $0\mapsto i\mapsto i-1\mapsto-i\mapsto i-1$, and
$x^p+c$ over $\QQ$ for $c\in\ZZ\setminus\{0\}$ when $p$ is odd, or $c\in\ZZ\setminus\{0,-1\}$ when
$p=2$.
\item every $(x+b)^2+c$ with $b,c\in\OK$ whose critical point $-b$ is not periodic; it is conjugate
to $x^2+(b+c)$ by $x\mapsto x+b$, and has the same $\kappa^*_\p$.
\end{enumerate}
\end{corollary}

\begin{proof}
(a), (b) follow from Lemma \ref{lem:kappaxd} and Example \ref{ex:kappa}. Over $\QQ$, $0$ is periodic for
$x^2+c$ only for $c\in\{0,-1\}$: if $f^{(q)}(0)=0$ the orbit of $0$ is bounded, so $-2\le c\le1/4$,
and $c=-2$ gives $0\mapsto-2\mapsto2\mapsto2$. For $x^p+c$ with $p$ odd and $c\in\ZZ\setminus\{0\}$:
if $|c|\ge2$ then $|w^p+c|\ge|w|^p-|c|>|w|$ whenever $|w|\ge|c|$, so $|f^{(n)}(0)|$ increases
strictly; if $c=\pm1$ the orbit is $0,c,2c,(2c)^p+c,\dots$ and the same argument applies from $2c$
on. (c) With $g(y)=f(y-b)+b=y^2+b+c$ we have $g'(y)=f'(y-b)$, and $y\mapsto y+b$ maps periodic points
and residue discs of $f$ to those of $g$; heights, fields and ramification indices are unchanged.
\end{proof}

\subsection{Explicit constants for \texorpdfstring{$x^2\pm1$}{x\^{}2+-1}}\label{sec:explicit}

Let $K=\QQ$, $p=d=2$ and $f=x^2+c$ with $c\in\{1,-1\}$, so $n_0=0$ and $q=2$ (Example
\ref{ex:n0q}). By Lemma \ref{lem:size} we may use $R=\phi=\frac{1+\sqrt5}2$ in Lemma
\ref{lem:basis}; then
$\Lambda=\frac12\log N+(\log_2N+1)\log\phi\le1.195\log N+0.482$.

For $x^2\pm1$ the bound for iterates of constant degree holds with the following explicit constants.

\begin{proposition}\label{prop:explicit}
Let $f=x^2\pm1$, let $\beta\in\Qbar$ be an algebraic integer with infinite $f$-orbit,
$D'=[\QQ(\beta):\QQ]$, and $\ell\ge\max\{\log6,\log3D'\}$. Put $L=\lceil12\sqrt\ell\,\rceil$.
If $[\QQ(f^{(n)}(\beta)):\QQ]=D'$ for $0\le n\le L$, then
$$\hh_f(\beta)\ \ge\ \frac1{4D'}\exp\bigl(-12\log2\cdot\sqrt\ell\bigr).$$
\end{proposition}

\begin{proof}
The proof follows the lines of the proof of Proposition \ref{prop:window}. With $\sum_{n\le L}2^n<2^{L+1}$,
\begin{equation}\label{eq:key}
\hh_f(\beta)\ \ge\ \frac{S_2(L)\log2-(L+1)\Lambda}{(L+1)\,D'\,2^{L+1}}.
\end{equation}
Here $\ell\ge\log6$ gives $L\ge17$. We claim
\begin{equation}\label{eq:claim}
S_2(L)\log2\ \ge\ 2(L+1)\Lambda.
\end{equation}
By \eqref{eq:Sq} it suffices that $(L+2)(L-3)\ge\frac{24}{\log2}\Lambda$. For $L\ge17$,
$(L+2)(L-3)\ge\bigl(1-\frac1{17}-\frac6{289}\bigr)L^2\ge0.92\cdot144\,\ell\ge132\,\ell$.
On the other side $\log N=\log(L+1)+\log D'$ with $\log D'\le\ell-\log3$ and
$L+1\le12\sqrt\ell+2\le12\ell+2$. So
$$\frac{24}{\log2}\Lambda\ \le\ 34.63\bigl(1.195(\ell-1.098+\log(12\ell+2))+0.482\bigr)
\ \le\ 41.4\,\ell-28.7+41.4\log(12\ell+2).$$
It remains to check $90.6\,\ell+28.7\ge41.4\log(12\ell+2)$ for $\ell\ge\log6$. At
$\ell=\log6$ the two sides are $191$ and $131$, and the derivative of the left side,
$90.6$, exceeds that of the right side, which is less than $22$. This proves
\eqref{eq:claim}. Then \eqref{eq:key} gives $\hh_f(\beta)\ge\Lambda/(D'2^{L+1})$, and
$\Lambda\ge1$ since $N\ge L+1\ge18$. With $L+1\le12\sqrt\ell+2$ we get
$\hh_f(\beta)\ge\frac1{4D'}2^{-12\sqrt\ell}$.
\end{proof}

\begin{proof}[Proof of Theorem \ref{thm:B}]
Let $e=e_2(\alpha)\le D$ and $\ell=\log3D$. If $\alpha$ is not an algebraic integer, Lemma \ref{lem:nonintegral} gives
$\hh_f(\alpha)\ge\log2/D$, which is larger than claimed; so let $\alpha$ be integral. The orbit of $0$ is
$0,1,2,5,\dots$, so no finite critical point is periodic, and $\kappa^*_2=2$ by Example \ref{ex:kappa}; thus
$\nu_2=2$. By Theorem \ref{thm:droploc} there are no drops at $n\ge m=\lceil\log_2e\rceil+2$. Proposition
\ref{prop:explicit}, applied to $\beta=f^{(m)}(\alpha)$ with $D'=[\QQ(\beta):\QQ]\le D$ and this $\ell$, gives
$\hh_f(\beta)\ge\frac1{4D}\exp(-12\log2\sqrt\ell)\ge\frac1{4D}\exp(-8.4\sqrt\ell)$. As $2^m\le4\cdot2^{1+\log_2e}=8e$,
we get $\hh_f(\alpha)=2^{-m}\hh_f(\beta)\ge(32eD)^{-1}\exp(-8.4\sqrt\ell)$.
\end{proof}

For $x^2-1$, Theorems \ref{thm:droploc} and \ref{thm:attract} with explicit constants give the following bound.

\begin{proposition}\label{prop:explicitminus}
Let $f=x^2-1$, let $\alpha\in\Qbar$ be not preperiodic with $D=[\QQ(\alpha):\QQ]\ge2$, and put $e=e_2(\alpha)$ and
$E_D=\exp(-8.4\sqrt{\log3D})$. Then
$$\hh_f(\alpha)\ \ge\ \min\Bigl\{\frac{E_D}{16\,e\,D},\ \frac{\log2}{3072\,e}\Bigr\}.$$
\end{proposition}

\begin{proof}
If $\alpha$ is not integral, Lemma \ref{lem:nonintegral} gives $\hh_f(\alpha)\ge\log2/D$, which exceeds the bound; so let $\alpha$ be
integral. The superattracting cycle is $\{0,-1\}$, so $\Sfin=\{\bar0,\bar1\}$, $Q=M=m_f=2$ and $\omega_f=2$. A periodic point $\eta$ with
$\bar\eta\notin\Sfin$ is a unit, so $\vartheta(\eta)=v(2\eta)=1$ and $\nu_2=1$. Since $|\log^+|x^2-1|_v-2\log^+|x|_v|\le\log2$ at the
archimedean place and the left side vanishes at the others, $|h(f(x))-2h(x)|\le\log2$, and Tate's telescoping argument
\cite[Thm.~3.20]{Sil07} gives $|\hh_f-h|\le\log2/(d-1)$; so we may take $C_f=\log2$, and then $A=2/\log2$, $B=3\log2$ and $2AB=12$ in
Theorem \ref{thm:attract}.

Suppose first that some $K$-conjugate of $\alpha$ reduces outside $\Sfin$, and let $m=\lceil\log_2e\rceil+1$, so that $2^m<4e$. By
Theorem \ref{thm:droploc} there are no drops at $n\ge m$, and Proposition \ref{prop:explicit}, applied to $f^{(m)}(\alpha)$ with
$\ell=\log3D$, gives $\hh_f(\alpha)\ge E_D/(4D_m2^m)\ge E_D/(16eD)$. If all $K$-conjugates of $\alpha$ reduce into $\Sfin$, Theorem
\ref{thm:attract} gives $\hh_f(\alpha)\ge\frac{3\log2}{64}\cdot\frac1{144}\cdot\frac1e=\frac{\log2}{3072\,e}$.
\end{proof}

For comparison, \cite[Cor.~1.4]{HS21} gives $\hh_f(\alpha)>0.007/D^2$ for every non-preperiodic $\alpha$. The bound $\log2/(3072e)$, valid
when every conjugate of $\alpha$ reduces into $\Sfin$, is better for all $D\ge32$, as $e\le D$; the bound $E_D/(16eD)$ is better only
when $e<8.9\,D\,E_D$, for instance when $e\le D^{1-\varepsilon}$ and $D$ is large. For points with $e\asymp D$ and a conjugate outside
$\Sfin$, the improvement on \cite{HS21} comes from Theorem \ref{thm:all}, whose constants we have not made explicit.

\end{document}